\documentclass[a4paper, 12pt]{article}
\usepackage{enumerate, theorem}
\usepackage{amsmath, amsfonts, amssymb}
\usepackage[height=22.5cm, width=15cm]{geometry}
\usepackage[all]{xy}
\usepackage{mathrsfs, yfonts}
\usepackage[utf8]{inputenc}

\newcommand{\A}{\tilde{\mathcal{A}}}
\newcommand{\Ab}{\mathcal{A}}
\newcommand{\Abc}{\tilde{\mathcal{A}}_{conv.}}

\newcommand{\parag}[1]{\paragraph{\sc{#1.}}}

\newtheorem{thm}{Theorem}[subsection]
\newtheorem{defn}[thm]{Definition}
\newtheorem{cor}[thm]{Corollary}
\newtheorem{prop}[thm]{Proposition}
\newtheorem{lemma}[thm]{Lemma}

\begin{document}
\title{Generalized Brieskorn Modules I: \\ Convergent (a,b)-modules}

\author{Daniel Barlet\footnote{Barlet Daniel, Institut Elie Cartan UMR 7502  \newline
Universit\'e de Lorraine, CNRS,  et  Institut Universitaire de France, \newline
BP 239 - F - 54506 Vandoeuvre-l\`es-Nancy Cedex.France. \newline
e-mail : daniel.barlet@univ-lorraine.fr}.}

\date{Mai 2026}

\maketitle

\hfill {\it A la mémoire de mon ami Alan T. Huckleberry}

%\hfill {\it    Les  hommes  d\'erivent  pour   comprendre} \\
%$.\qquad \qquad \qquad  \qquad\qquad \qquad \qquad\qquad \qquad  \quad$ {\it  alors  que  Dieu se  contente  d'int\'egrer}.

\parag{Abstract} This paper is the first one of two papers whose  goal is to give a converse to the main result of my  previous paper \cite{[B.22]}, so to prove the existence of multiple poles for the distribution $\vert f\vert^{2\lambda}$ with an hypothesis on a {\bf Higher  Bernstein Polynomial} of the (a,b)-module  generated by  the germ $\omega \in \Omega^{n+1}_0$ of a given holomorphic volum form. Note that, even for the existence of  a simple pole this converse is already new.  One difficulty to prove such a result comes from  the use of the formal completion in $f$ of the Brieskorn module of the holomorphic germ $f: (\mathbb{C}^{n+1}, 0) \to (\mathbb{C}, 0)$ which  does not give access to the cohomology of the Milnor's fiber of $f$, which  by definition,  is outside $\{f=0\}$. This leads to introduce  {\bf generalized Brieskorn modules} (convergent geometric  (a,b)-modules)  which allow this passage. The first aim of this part I is to give a solid basis of the theory of convergent (a,b)-modules. \\
In order to take in account Jordan blocs of the monodromy in our results we introduce the {\bf semi-simple filtration} of a {\bf generalized Brieskorn module} (convergent  (a,b)-module) and we shall use it to define in part II   the {\bf higher order Bernstein polynomials} in this context. They correspond to a decomposition of the ``standard'' Bernstein polynomial  of a generalized Brieskorn module, taking in account the nilpotent order of the monodromy. \\
In this part I we obtain also a  full description of generalized Brieskorn-modules in terms of (convergent) asymptotics expansions of Nilsson class which will be used as a starting point in part II. We conclude this part  I by making explicite  the relationship between the semi-simple filtration of a generalized Brieskorn module $E$  and the nilpotent filtration  of the monodromy on its saturation $E^\sharp$.

\newpage
\parag{AMS classification}  32 S 25; 32 S 40 ; 34 E 05

\parag{Key words} Generalized Brieskorn Module;  Bernstein Polynomial; Asymptotic Expansions ; Period-integral;  Convergent (a,b)-Module; Geometric  (a,b)-Modules;  Semi-Simple Filtration;  Fresco. \\

\tableofcontents

\newpage

\section{Introduction}

\subsection{Preamble}

After having worked during the 80's on the singularities of hyper-surfaces and the associated  regular meromorphic connections (Brieskorn modules), I wanted to understand what was the results on these singularities which are consequences only of the  action of "a" (product by $f$) and "b" (the inverse of derivation by $f$). Indeed in his article \cite{[K-S]} Kyoji Saito remarks that $\partial^{-1}$ acts on the Brieskorn lattice of an isolated singularty of hypersurface. This leads me to study the  non commutative algebraic structure that I called (a, b)-module. Since the regularity the the Gauss-Manin connection implies that any formal solution is in fact convergent, I decide to adopt a formal context (in "b", but which allows also formal series in "a") for regular (a, b)-modules.\\
Then I define the notion of Bernstein polynomial for a regular (a,b)-module. M. Kashiwara's Theorem which shows that the roots of the Bernstein polynomial of a holomorphic function are always negative rational numbers justifies the introduction of "geometric (a, b)-modules" which are simply regular (a,b)-module with this property.\\
Although this setting allows me to obtain several new results with this approach (see for instance \cite{[B.93]}, \cite{[B.II]},  \cite{[B.09]}, \cite{[B.22]}) for more than twenty years, the results I want to prove these last years and that are obtained in \cite{[GBM-II]} force me to introduce a convergent version of the geometric (a,b)-modules that I call Generalized Brieskorn Modules. \\
So the first goal of this article  is to introduce convergent series (in fact Gevrey convergent series in "b") to be able to complete my proofs. \\
In  the begining  of this article I  give the estimates useful to obtain the basic results in the convergent frame for (a,b)-modules including in Section 5  the relationship with convergent  asymptotics expansions.\\
In Section 4  is  established the existence of the semi-simple filtration of a generalized Brieskorn module which will be used  in \cite{[GBM-II]} to give a decomposition of its Bernstein polynomials in { \bf higher order Bernstein polynomials} which encode, for each root  $r $ of the Bernstein polynomial, its effect on nilpotent order  of the corresponding eigenvalue $exp(-2i\pi r)$ of the  monodromy associated  to the generalized Brieskorn module. \\
The relationship between the semi-simple filtration is described in Section 6.  \\

Note that in the present article and in the second part \cite{[GBM-II]} we use mainly the action on generalized Brieskorn modules of  the algebra $B[a]$ of polynomials in $"a"$ with coefficient in the algebra $B := \mathbb{C}\{\{b\}\}$ (Gevrey series in $"b"$) with the commutation relation $aS(b) = S(b)a + b^2S'(b)$ where $S'$ is the (usual) derivative of $S$. This is enough for the results in \cite{[GBM-II]}. But this is not completely satisfying because we have not only  a natural action of $\mathbb{C}[a]$ on a generalized Brieskorn module, but also a natural action of $\mathbb{C}\{a\}$. This leads to defined the "good" algebra acting on generalized Brieskorn modules and this is the goal of the article \cite{[GBM-III]}

\subsection{Descrition of the content}

 \begin{itemize}
\item In section 2 we establish several  basic estimates to work with the algebra $B := \mathbb{C}\{\{b\}\}$ which is a closed  sub-algebra of  the algebra of continuous $\mathbb{C}$-linear endomorphisms  $End_c(\mathbb{C}\{s\})$ of the algebra $\mathbb{C}\{s\}$ of holomorphic germs at the origin in $\mathbb{C}$.\\
We shall denote by $A$ the sub-algebra of $End_c(\mathbb{C}\{s\})$ given by multiplication by an element in $\mathbb{C}\{s\}$. Denoting by $a$ the multiplication by $s$ we have $A = \mathbb{C}\{a\}$.\\
We conclude this section by a construction which an analog "(a,b)-module construction" of the Deligne's construction in \cite{[D.70]}.
 \item Section 3 is devoted to extend the standard properties of regular (formal)  (a,b)-modules (see \cite{[B.93]}) to the convergent case. 
 \item The purpose of Section 4 is to define the semi-simple filtration of a regular convergent (a,b)-module and to establish the basic properties of this filtration. 
 \item In Section 5 we show that convergent geometric  (a,b)-modules are simply convergent sub-(a,b)-modules of classical  convergent asymptotic expansions modules. The Embedding Theorem shows that the geometric (a,b)-modules are exactly  regular (a,b)-modules such that the roots of their Bernstein polynomial are negative and rational (compare with with the positivity Theorem of B. Malgrange \cite{[M.74]} and the Rationality Theorem of M. Kashiwara \cite{[K.76]}). 
  \item Section 6  explains the relation between the semi-simple filtration on a regular (convergent) (a,b)-module $\mathcal{E}$ and  the nilpotent part of the  monodromy  which is naturally defined on a simple pole geometric (a,b)-module, so on the saturation $\mathcal{E}^\sharp$ of $\mathcal{E}$. 
\end{itemize}

\section{Convergent (a,b)-modules}

\parag{Notations} In the sequel we shall use the following complex unitary algebras:
\begin{enumerate}
\item We shall denote by $A$ the sub-algebra of $End_c(\mathbb{C}\{s\})$ given by multiplication by an element in $\mathbb{C}\{s\}$. Denoting by $a$ the multiplication by $s$ we have $A = \mathbb{C}\{a\}$.
\item The algebra $B := \mathbb{C}\{\{b\}\}$ which is the closed  sub-algebra of  continuous linear operators on $\mathbb{C}\{s\}$  generated by $b$, the primitive without constant in $\mathbb{C}\{s\}$.
\item The algebra $\mathcal{A}$ which is the algebra of polynomials in $a$ and $b$ inside the algebra of continuous linear operators on $\mathbb{C}\{s\}$. So they satisfy  $ab - ba = b^2$.
\item The algebra $B[a]$ of polynomials in $a$ with coefficients in $B$, with the commutation relations $aS(b) - S(b)a = b^2S'(b)$ for each $S \in B$.
\item The algebra $\widehat{A} := \mathbb{C}[[a]]$
\item The algebra $\widehat{B} := \mathbb{C}[[b]]$
\item The algebra $\tilde{\mathcal{A}} := \widehat{B}[a]$ of polynomials in $a$ with coefficients in $\widehat{B}$, with the commutation relations $aS(b) - S(b)a = b^2S'(b)$ for each $S \in \widehat{B}$.
\item The algebra $\widehat{\mathcal{A}}$ which the algebra of formal power series in $a$ and $b$ with the commutation relation $ab - ba = b^2$.
\end{enumerate}

\subsection{The algebra $\mathbb{C}_r\{\{b\}\}$}

For each real number $r \in ]0, 1[$ let $\mathbb{C}_r\{\{b\}\}$ the sub-vector space of $\mathbb{C}[[b]]$ defined by
\begin{equation}
\mathbb{C}_r\{\{b\}\} := \Big\{ S = \sum_{j=0}^\infty  s_jb^j \ / \ \exists R \in ]1,1/r[, \ \exists C_R \ :  \quad   \vert s_j\vert \leq C_RR^jj!  \quad \forall j \in \mathbb{N}  \Big\}.
\end{equation}

\parag{Remarks}\begin{enumerate}
\item If the condition above holds for some $R_0 \in]1,1/r[$ then it holds for any $R$ in $[R_0, 1/r[$ and with the constant $C_R = C_{R_0}$.
\item If for some positive integer  $q $ and  for some $ S = \sum_{j=0}^\infty  s_jb^j$ we have:
$$ \ \exists R \in ]1,1/r[ \quad \exists C_R \ {\rm such \ that} \quad \vert s_j\vert \leq C_RR^j(j+q)! \quad \forall j \in \mathbb{N}$$
then $S$ is in $\mathbb{C}_r\{\{b\}\}$ thanks to the easy estimates
 $$\frac{(j+q)!}{j!} \leq (1+q)^q j^q \quad \forall q \in \mathbb{N}, \forall j \in \mathbb{N}^*,$$
the sequence $j \mapsto (R/\rho)^jj^q$ being bounded for $q$ fixed and  for  $\rho \in ]R, 1/r[$.
 \item The vector space $\mathbb{C}_r\{\{b\}\} $ is stable by derivation because if $S = \sum_{j=0}^\infty  s_jb^j$ is in $\mathbb{C}_r\{\{b\}\} $  then $S' := \sum_{j=0}^\infty (j+1)s_{j+1}b^j$ satisfies
 $$\exists R \in ]1,1/r [ \ : \quad  \vert  (j +1)s_{j+1} \vert \leq  (j+1)C_R R^{j+1} (j+1)! \leq RC_R R^j (j+2)! $$
  and   we may apply the previous remark with $q = 2$ to $S'$.$\hfill \square$\\
\end{enumerate}

We have a vector spaces isomorphism:
$$\mathbb{C}_r\{\{b\}\} \to \mathbb{C}_r\{z\} := \Big\{ \sum_{j=0}^\infty  c_jz^j \ / \  \exists R \in ]1, 1/r[ , \ \exists C_R \ :   \quad \ \vert c_j\vert \leq C_RR^j \quad \forall j \in \mathbb{N} \Big\}$$
which is given by $s_j \mapsto c_j = s_j/j!$. As $\mathbb{C}_r\{z\}$ is the algebra of germs of holomorphic functions around the closed disc $\bar D_{r}$ in $\mathbb{C}$, it is a {\bf dual Fréchet} (in short DF)  topological  vector space and we shall define the topology on $\mathbb{C}_r\{\{b\}\}$ via this isomorphism.

\begin{lemma}\label{conv. 1}
For each $r \in ]0, 1[$, $\mathbb{C}_r\{\{b\}\}$ is a dual Fr\'echet  sub-algebra of $\mathbb{C}[[b]]$ which is stable by derivation.
\end{lemma}

\parag{proof} Let $S :=  \sum_{j=0}^\infty  s_jb^j $ and $T :=  \sum_{j=0}^\infty  t_jb^j $ be in  $\mathbb{C}_r\{\{b\}\}$. The product $ST$ in $\mathbb{C}[[b]]$ is given by
$$ ST = \sum_{j=0}^\infty  u_jb^j \quad {\rm where} \quad  u_j := \sum_{p=0}^j  s_{j-p}t_p .$$
For some  $R \in ]1, 1/r[$ there exists $C_R > 0$ and $D_R > 0$ such that $\vert s_j\vert \leq C_RR^jj!$ and $\vert t_j\vert \leq D_RR^jj!$ by definition. This implies
$$ \vert u_j\vert \leq C_RD_RR^jj!\big(\sum_{p=0}^j \frac{(j-p)!p!}{j!}\big).$$
The following easy lemma will be used several times in the sequel.

\begin{lemma}\label{inegalite 1}
 For any $j \geq 0$ we have \  $\sum_{p=0}^j (j-p)!p! \leq 3(j!)$.
 \end{lemma}
 
\parag{Proof}  For $j \geq 3$ and $p \in [1, j-1]$ we have   $(j-p)!p! \leq (j-1)!$, so  for $j \geq 3$:
$$ \sum_{p=0}^j (j-p)!p! \leq 2(j!) +  \sum_{p=1}^{j-1} (j-p)!p! \leq 2(j!) + (j-1)((j-1)!) \leq 3(j!) .$$
As the cases $j = 0, 1, 2$ are obvious, so the proof is complete. $\hfill \blacksquare$

\parag{End of proof of Lemma \ref{conv. 1}} So we obtain that $ \vert u_j\vert \leq 3C_RD_RR^jj!$ showing that $ST$ is in $\mathbb{C}_r\{\{b\}\}$.\\
The countable family given by  
 $$\vert\vert S\vert\vert_R := \sup_j \vert s_j\vert\big/(R^jj!)$$ 
  for $R \in \mathbb{Q}\cap ]1, 1/r[$ defines the dual Fr\'echet topology on $\mathbb{C}_r\{\{b\}\}$ and the previous computation
shows that $\vert\vert ST\vert\vert_R \leq 3\vert\vert S\vert\vert_R\vert\vert T\vert\vert_R$ for each $R ]1,1/r[$.\\
The stability by derivation is explained in Remark 3 above.$\hfill \blacksquare$\\

For each real number  $r \in [0, 1[$ define  the operator $b$ on the algebra  $\mathbb{C}_r\{z\}$  of  germs of holomorphic functions around the closed disc $\bar D_{r}$  with center $0$ and radius $r$ in $\mathbb{C}$, as the linear continuous operator given by the primitive vanishing at the origin;  that is to say:
  $$b(f)(z) := \int_0^z f(t)dt = z\int_0^1 f(sz)ds.$$

\begin{prop}\label{conv. 2}
This action of $\mathbb{C}[b]$ on $\mathbb{C}_r\{z\}$ extends to a continuous action of $\mathbb{C}_r\{\{b\}\}$ on $\mathbb{C}_r\{z\}$. It induces an isomorphism of  DF-algebras between  $\mathbb{C}_r\{\{b\}\}$ and  the commutant of $b$ in the algebra of continuous endomorphism of $\mathbb{C}_r\{z\}$.
\end{prop}

Note that this commutant is a closed sub-algebra of the algebra of continuous endomorphism of $\mathbb{C}_r\{z\}$.

\parag{proof}  Let $S := \sum_{j=0}^\infty  s_jb^j$ and $f(z) := \sum_{\nu = 0}^\infty  \gamma_\nu z^\nu$ where we may assume that there exists $R \in ]1,1/r[$ such  that  the following estimates hold:
  $$ \vert s_j\vert \leq C_RR^jj! \quad {\rm and}\quad  \vert \gamma_\nu \vert \leq \Gamma_R R^\nu .$$
  Since  $b^j[z^\nu] = \frac{\nu! z^{\nu+j}}{(\nu+j)!}$ we obtain, writing  $S[f](z) = \sum_{q=0}^\infty \delta_qz^q$
  $$\vert \delta_q \vert =  \vert \sum_{j=0}^q s_j \gamma_{q-j}\frac{(q-j)!}{q!} \vert \leq \frac{C_R\Gamma_R R^q}{q!} \sum_{j=0}^q  j!(q-j)! .$$
 We obtain, thanks to the estimates given in Lemma \ref{inegalite 1}
  $$ \vert \delta_q \vert \leq 3C_R\Gamma_RR^q$$
  This shows that $S[f]$ is in $\mathbb{C}_r\{z\}$ and that $S$ acts continuously on this DF space. \\
  Then it is clear that the image of $\mathbb{C}_r\{\{b\}\}$  is a commutative sub-algebra of the commutant of $b$ in the algebra of  continuous endomorphisms of $\mathbb{C}_r\{z\}$.\\
  Let us show the converse.
  If $\Phi$ is a continuous endomorphisms of $\mathbb{C}_r\{z\}$ commuting with $b$, put $\Phi(1) = \sum_{j=0}^\infty \gamma_jz^j$ and define $S := \sum_{j=0}^\infty j!\gamma_jb^j$. Then $S$ is in
   $\mathbb{C}_r\{\{b\}\}$. We shall compare $\Phi[z^p]$ and $S[z^p]$ for $p \in \mathbb{N}$:
  $$ \Phi[z^p] = \Phi[p!b^p(1)] = p!b^p[\Phi(1)] = \sum_{j=0}^\infty p!\gamma_j\frac{j!}{(p+j)!}z^{p+j} $$
  and
  $$ S[z^p] = \sum_{j=0}^\infty j!\gamma_jb^j(z^p) =  \sum_{j=0}^\infty j!\gamma_j\frac{p! z^{p+j}}{(p+j)!}\ \cdot$$
  As the linear combinations of the $z^p, \, p \in \mathbb{N}$, are dense in $\mathbb{C}_r\{z\}$ we conclude that $\Phi$ coincides with the image of $S$ in the algebra of  continuous endomorphisms of $\mathbb{C}_r\{z\}$. So the image of $\mathbb{C}_r\{\{b\}\}$ is the commutant sub-algebra of $b$.$\hfill \blacksquare$\\
  
   Note that the action of $\mathbb{C}_r\{\{b\}\}$ extends to the Banach algebra of  continuous functions on $\bar D_r$ using the formula
    $$b(f)(z) :=  z\int_0^1 f(tz)dt $$
    since it  is easy to see that for a continuous function  on $\bar D_r$ we have the estimate
     $$\vert\vert b^j(f)\vert\vert_{r} \leq \vert\vert f\vert\vert r^j/j! .$$

  For each real number $s \in ]0, r[$ we have a continuous inclusion  $\mathbb{C}_s\{\{b\}\} \subset \mathbb{C}_r\{\{b\}\}$ of DF-algebras and we define $B$ as the algebra $\mathbb{C}_0\{\{b\}\}$ which is the union for all positive $r$ of the algebras $\mathbb{C}_r\{\{b\}\} $. \\
  Note that the algebra $B$ is defined by
 $$B := \Big\{S = \sum_{j=0}^\infty  s_jb^j \ / \ \exists R > 1, \   \exists  C_R > 0 \ : \quad  \vert s_j\vert \leq C_RR^j j!  \quad \forall j \in \mathbb{N} \Big\}$$
 and that $B$ acts on $ \mathbb{C}_0\{z\}$, the algebra of holomorphic germs at the origin which is the union for all $r \in ]0, 1[$ of the algebra $\mathbb{C}_r\{z\}$ and is defined by
 $$ \mathbb{C}\{s\} := \Big\{f = \sum_{j=0}^\infty  \gamma_jz^j \ / \ \exists R > 1, \  \exists  C_R > 0 \ :  \quad  \vert \gamma_j\vert \leq C_RR^j \quad \forall j \in \mathbb{N}  \Big\}.$$
  
  \begin{prop}\label{conv.3bis}
  The algebra $B$ is local.
  \end{prop} 
  
 \parag{Proof}  This result is a consequence of the following lemma.$\hfill \blacksquare$\\
  
  \begin{lemma} \label{conv. 3}
  Let $S$ be in $\mathbb{C}_r\{\{b\}\}$ such that $S(0) \not= 0$. Then there exists $r' \in ]0,r[$ such that $S$ is invertible in $\mathbb{C}_{r'}\{\{b\}\}$. So the algebra $B := \mathbb{C}_0\{\{b\}\} $ has a unique closed maximal ideal which is  generated by $b$.
    \end{lemma}
  
  \parag{proof} The inverse of $S$ in $\mathbb{C}[[b]]$ is given by $T := \sum_{j=0}^\infty t_jb^j$ with $t_0 := 1/s_0$ and $t_q := -(1/s_0)(\sum_{j=1}^q s_jt_{q-j})$ for $q \geq 1$. For some   $R \in ]1, 1/r[$ there exists $C_R > 0$ such that $\vert s_j\vert \leq C_RR^jj! \quad \forall j \geq 0$. Now assume that for some \\
   $$\rho > \sup\big\{ \frac{1}{r}\ , \frac{3RC_R}{\vert s_0\vert} \big\}, $$
some $D_\rho > 0$ and some $q \in \mathbb{N}^* $    the estimates  $ \vert t_h\vert \leq D_{\rho}\rho^hh! $  holds true for each integer  $h \in [0, q-1]$.
 Then we obtain
   \begin{align*}
   & \vert t_q\vert \leq \frac{1}{\vert s_0\vert}C_RD_{\rho}\sum_{j=1}^q R^j\rho^{q-j}j!(q-j)! \\
   &  \vert t_q\vert \leq \frac{1}{\vert s_0\vert}C_RD_{\rho}\rho^qq! \sum_{j=1}^q (\frac{R}{\rho})^j\frac{j!(q-j)!}{q!} \\
   &  \vert t_q\vert \leq \frac{1}{\vert s_0\vert}C_RD_{\rho}\rho^qq! \frac{3R}{\rho} \leq D_{\rho}\rho^qq!
   \end{align*}
   So the estimates $\vert t_h\vert \leq D_{\rho}\rho^hh! $ will be valid $\forall h \geq 0$  as soon as it is true for $h = 0$. We conclude the proof by defining $D_\rho := 1/\vert s_0\vert$ and by choosing
    $r' < 1/\rho $.$\hfill \blacksquare$ \\
   
   \subsection{The r-convergent (a,b)-modules}

   Fix a real number $r \in [0, 1[$. 
   
   \begin{defn}\label{cv 0}
   A free finite rank \  $\mathbb{C}_r\{\{b\}\}$-module $\mathcal{E}_r$ endowed with a continuous $\mathbb{C}$-linear endomorphism $a$ which satisfies
   $$ ab - ba = b^2$$
   will be called a {\bf r-convergent (a,b)-module}. In the case $r = 0$ we simply call $\mathcal{E} := \mathcal{E}_0$  a {\bf convergent (a,b)-module}.\\
   We say that $\mathcal{E}_r$  has  a {\bf simple pole}  when it satisfies $a\mathcal{E}_r \subset b\mathcal{E}_r$.
   \end{defn}
   
   \parag{Remark} The dual Fr\'echet topology on a free finite rank $\mathbb{C}_r\{\{b\}\}$-module $\mathcal{E}_r$   given by the choice of a $\mathbb{C}_r\{\{b\}\}$-basis  $e$ of $\mathcal{E}_r$   is independent of this choice because if $\varepsilon := M(b)e$ is an other basis, the linear bijective map corresponding to the change of basis is continuous, thanks to the continuity of the product in the $\mathbb{C}$-algebra $\mathbb{C}_r\{\{b\}\}$. Then it is an isomorphism of dual Fr\'echet spaces.\\
   
   Note that the continuity of $a$ for the natural dual Fr\'echet topology deduced from any $\mathbb{C}_r\{\{b\}\}$-basis of $\mathcal{E}_r$  implies that for any $S(b) := \sum_{j=0}^{+\infty} s_jb^j $ \ in $  \mathbb{C}_r\{\{b\}\}$ and for any $x \in \mathcal{E}_r$  we have
  $$a(S(b)x) = \lim_{N \to \infty} \sum_{j=0}^N s_ja(b^jx) = \lim_{N \to \infty}  \sum_{j=0}^N s_j(b^ja + jb^{j+1})x = S(b)ax + b^2S'(b)x .$$
  Of course, to each r-convergent (a,b)-module we can associate a  s-convergent (a,b)-module for any $s \in [0, r]$ and so a convergent (a,b)-module in the case $s = 0$, via the correspondence:
  $$\mathcal{E}_r \mapsto \mathcal{E}_r \otimes_{\mathbb{C}_r\{\{b\}\}} \mathbb{C}_s\{\{b\}\} $$
where  the action of $a$ is  defined by
 $$a(x\otimes S) = ax\otimes S + x\otimes b^2S'$$
 for $x \in \mathcal{E}_r$.
  
     \begin{lemma}\label{cv 2}
   For any convergent (a,b)-module with a simple pole $\mathcal{E}$  there exists  $r > 0$  and a r-convergent  (a,b)-module $\mathcal{E}_r \subset \mathcal{E}$ such that
    $\mathcal{E} = \mathcal{E}_r \otimes_{\mathbb{C}_r\{\{b\}\}} B$ as a   $B$-module and such that the equality
   $$ a(x\otimes S(b)) = ax\otimes S(b) + x\otimes b^2S'(b) $$
   holds for each $S \in B$ and each $x \in \mathcal{E}_r$.
   \end{lemma}
   
   \parag{proof}  Let $e := (e_1, \dots, e_k)$ be a  $B$-basis of $\mathcal{E}$ and write $ae = M(b)e$ where $M$ is in $B\otimes_{\mathbb{C}} End_{\mathbb{C}}(\mathbb{C}^k)$ satisfies 
   $M(0) = 0$. Choose now $r > 0$ such that $M$ is in fact in $\mathbb{C}_r\{\{b\}\}\otimes_{\mathbb{C}} End_{\mathbb{C}}(\mathbb{C}^k)$ and define $\mathcal{E}_r := \oplus_{j=1}^k \mathbb{C}_r\{\{b\}\}e_j \subset \mathcal{E}$. Now the $\mathbb{C}$-linear endomorphism of $\mathcal{E}_r$ induced by $a$ is clearly continuous for the dual Fr\'echet topology of $\mathcal{E}_r$, thanks to the continuity of the product and of the derivation in $\mathbb{C}_r\{\{b\}\}$, and satisfies $ab - ba = b^2$. It is easy to verify that the r-convergent  (a,b)-module $\mathcal{E}_r$ satisfies the lemma.$\hfill \blacksquare$

 \subsection{Construction of $\mathcal{E}(\Theta)$} 
 
 We shall construct now an important family of examples of simple poles convergent (a,b)-modules.\\
 We begin by  a  very simple but  useful lemma:

\begin{lemma}\label{simple}
For $x \in [1, +\infty[$ and any integer $k \geq 0$ we have
$$ \frac{x(x+1)\dots (x+k)}{(k+1)!} \leq x^{k+1}.$$
\end{lemma}

\parag{Proof} Obvious because for $x \geq 1$ we have $x+j \leq x(j+1)$ for each $j \geq 0$.$\hfill \blacksquare$\\

Note that for $x \in [0, 1]$ we have 
$ \frac{x(x+1)\dots (x+k)}{(k+1)!} \leq  1.$\\

As an easy consequence, we obtain that for any endomorphism $\Theta \in L(\mathbb{C}^p, \mathbb{C}^p)$ we have the estimates
$$ \vert\vert \Theta\circ(\Theta+1)\circ \dots\circ (\Theta+ k) \vert\vert \leq (k+1)!\,\theta^{k+1} $$
where $\theta := sup \{1, \vert\vert \Theta\vert\vert \}$ and where the norm $\vert\vert\quad  \vert\vert$ satisfies $\vert\vert x\circ y\vert\vert \leq \vert\vert x\vert\vert.\vert\vert y \vert\vert $ for any $x, y \in L(\mathbb{C}^p, \mathbb{C}^p)$.\\

Then if we have a basis $e$ of a rank $p$ convergent  (a,b)-module which satisfies
$$ ae = \Theta be $$
where $\Theta$ is in $L(\mathbb{C}^p, \mathbb{C}^p)$ we shall have, by an easy induction on $k \geq 1$:
$$ a^ke = \Theta\circ(\Theta+1)\circ\dots \circ (\Theta+k)b^ke $$
because $ab^k = b^ka + k b^{k+1}$. \\
Assume now that we have $T \in \mathbb{C}\{a\}\otimes_{\mathbb{C}} L(\mathbb{C}^p, \mathbb{C}^p)$ and $x = Te$;   write \\
 $T := \sum_{k=0}^\infty T_ka^k$ with a positive radius of convergence. We obtain :
$$ x = \sum_{k=0}^\infty T_k\Theta\circ(\Theta+1)\circ\dots\circ(\Theta+k)b^ke $$
and the estimates above implies that  $x$ lies inside the $B$-module generated by $e$, because the series $\frac{1}{k!}T_k\Theta\circ(\Theta+1)\circ\dots\circ(\Theta+k)b^k$ has a positive radius of convergence.\\

\begin{prop}\label{cv 3}
Let $\Theta$ be an invertible $(k,k)-$matrix with complex entries such that the spectrum of $\Theta$ is disjoint from $-\mathbb{N}$. Then we define the convergent (a,b)-module  $\mathcal{E}(\Theta)$  as follows:
\begin{enumerate}
\item $\mathcal{E}(\Theta)$ is the free, rank $k$, $B$-module with basis $e := (e_1, \dots, e_k)$, the standard basis of $\mathbb{C}^k$.
\item The $\mathbb{C}$-linear {\bf continuous} endomorphism $a : \mathcal{E}(\Theta) \to \mathcal{E}(\Theta)$ is defined by 
$$ ae = \Theta be \quad {\rm and} \quad ab - ba = b^2.$$
\end{enumerate}
Then $\mathcal{E}$ is a simple pole convergent (a,b)-module. \\
Moreover, the action of $\mathbb{C}[a]$ on $\mathcal{E}(\Theta)$ extends to a continuous action of $A = \mathbb{C}\{a\}$ on $\mathcal{E}(\Theta)$ and under this action $\mathcal{E}(\Theta)$ is a free, rank $k$, $A$-module with basis $e$.
\end{prop}

\parag{proof} In fact we shall begin by proving that for any $r \in ]0, 1[$ the $\mathbb{C}_r\{\{b\}\}-$sub-module $\mathcal{E}_r(\Theta) := \oplus_{j=1}^k \mathbb{C}_r\{\{b\}\}e_j \subset \mathcal{E}(\Theta)$ which is clearly stable by $a$, is  a r-convergent (a,b)-module with a simple pole. The continuity of $a$  for the dual Fr\'echet topology of $\mathcal{E}_r$ follows from the continuity of the derivation in $\mathbb{C}_r\{\{b\}\}$ and the formula  $a(S(b)e) = S(b)\Theta e + b^2S'(b)e$ for any matrix 
$S \in \mathbb{C}_r\{\{b\}\}\otimes_{\mathbb{C}}End_{\mathbb{C}}(\mathbb{C}^k)$  .\\
The last assertion is an easy consequence of the following equality, valid for all $j $ and $p$ in $ \mathbb{N}$:
\begin{equation}
a^jb^pe = (\Theta+pId)\circ \dots\circ (\Theta + (p+j-1)Id)b^{(j+p)}e .
\end{equation}
As this equality is clear for $j = 0, 1$ and any $p \in \mathbb{N}$, thanks to the commutation  relation $ab^p = b^pa + pb^{p+1}$, assume that it is true for $j$ and compute $a^{j+1}b^pe$. The matrix $\Theta$ has complex entries so commutes with the action of $a$, so the commutation relation above is enough to conclude.\\
Now if $U := \sum_{j=0}^{+\infty} u_ja^j $ is in $\mathbb{C}_r\{a\}$, for some $R \in ]1, 1/r[$ we may find a constant $C_R > 0$ such that $\vert u_j\vert \leq C_RR^j \quad \forall j \geq 0$. 

Also if $S := \sum_{j=0}^{+\infty} s_jb^j$ is in $ \mathbb{C}_r\{\{b\}\}^k$ we may find, choosing a bigger $R \in ]1, 1/r[$ if necessary, a constant $D_R$ such that
 $$\vert\vert s_j\vert\vert \leq D_RR^jj!.$$
 Then we have
$$ U(S(b)e) = \sum_{p=0}^{+\infty}  X_pb^pe $$
where $X_p $ in $\mathbb{C}^k$ is given by
$$ X_p = \sum_{j=0}^p s_ju_{p-j}.$$
So we obtain the estimates, for some $R \in ]1, 1/r[$, where $\theta := sup \{1, \vert\vert \Theta\vert\vert\}$,  assuming that $\vert\vert \ \vert\vert$ is a multiplicative norm on the  complex $(k, k)$-matrices and  that $q$ is a positive integer bigger than $\theta$
$$ \vert\vert X_p\vert\vert \leq C_RD_RR^p\big(\sum_{j=0}^p \frac{(p+q-1))!}{(q+j-1)!}\big) \leq (1/(q-1)!) C_RD_RR^p(p+q)!$$
since we have
 $$\sum_{j=0}^p (q+p-1)!/(q+j-1)! \leq (p+q)!/(q-1)!$$
 for each integer $p \geq 0$. So $U(S(b)e)$ is in $\mathcal{E}_r(\Theta)$. This implies that $\mathcal{E}_r(\Theta)$ is a $\mathbb{C}_r\{a\}$-module.\\
 To show that $\mathcal{E}(\Theta)$ is free $A$-module with basis $e$, since  it has no $a$-torsion thanks to our hypothesis that $Spec(\Theta) \cap \{ -\mathbb{N}\} = \emptyset$, it is enough to show that, for  any  $S \in B\otimes_{\mathbb{C}} End_{\mathbb{C}}(\mathbb{C}^k)$, $S(b)e$ is in the $A$-module generated by $e_1, \dots, e_k$. Using the formula above with $p= 0$ we obtain, with the notation $(\Theta + qId)^{-1} := H_q$
$$ b^je = \big(H_{j-1}\circ \dots\circ H_{0}\big)a^je $$
for each $j \geq 0$ and then $S(b)e = \sum_{p=0}^{+\infty}  Y_pa^pe $ with $Y_p := s_p\big(H_{p-1}\circ \dots\circ H_{0}\big)$.\\
But since the sequence $(\Theta/q)$ converges to $0$ when $q \to + \infty$, there for each choice of $\varepsilon > 0$  exists a constant $\Gamma \in ]1, 1 + \varepsilon[ $ such that $(\Theta + qId)^{-1} = (1/q)(1 + \Theta/q)^{-1}$ has its norm bounded by $\Gamma/(q+1)$ for each $q$ large enough (but depending on the choice of $\varepsilon$). This implies $$\vert\vert H_{p-1}\circ \dots\circ H_{0}\vert\vert \leq \Delta(\varepsilon)\Gamma^p/p! \quad \forall p \geq 0 $$
 where the positive  constant $\Delta(\varepsilon)$ depends on the choice of $\varepsilon$. Since for some $R$  in $]1, 1/r[$ there exists a constant $C_R$ such that $\vert\vert s_p\vert\vert \leq C_RR^pp! \quad \forall p \geq 0$, we obtain the estimates:
$$\vert Y_p\vert \leq \Delta(\varepsilon) C_RR^p\Gamma^p $$
Now choose $\varepsilon$ such that $R(1+\varepsilon) < 1/r$ so that $R' := R\Gamma $ is in $]1, 1/r[$. Then the estimates above with $C_{R'} := \Delta(\varepsilon)C_R$
 implies that $\sum_{p=0}^{+\infty} Y_pa^p$ is in $\mathbb{C}_{r}\{a\} \subset A$.$\hfill \blacksquare$\\

\parag{Example} For $\alpha \in ]0,1] \cap \mathbb{Q}$, the convergent (a,b)-module $\Xi_\alpha^{(N)}$ of asymptotic expansions (see section 5  below)   is equal to  $\mathcal{E}(\Theta_\alpha)$ where the matrix $\Theta_\alpha$ is the following $(N+1,N+1)$ matrix:
$$J_{\alpha, N+1} := \begin{pmatrix} \alpha & 0 & 0 & 0 & \dots & 0 \\
 1 & \alpha & 0 & 0 &\dots &0 \\
  0   & 1& \alpha & 0 & 0  & 0 \\
  \dots &\dots & \dots &\dots &\dots &0 \\ 
0 & \dots &\dots & 1 &\alpha &0  \\
 0 & \dots & 0 & 0 & 1 & \alpha \end{pmatrix}.$$

\section{Regular convergent (a,b)-modules}

\subsection{Basic properties}

\parag{Notation} Recall that we  note $B$ and $\widehat{B}$ respectively the algebras $B := \mathbb{C}\{\{b\}\}$ and $\widehat{B} := \mathbb{C}[[b]]$ in the sequel.\\
For $S$ in $B$ (or in $\widehat{B}$) we shall note $S'$ the usual derivative of $S$; we mean that  $S'(b) = \sum_{j=1}^\infty js_jb^{j-1}$ if $S := \sum_{j=0}^\infty s_jb^j$. The algebras $B$ and 
$\widehat{B}$ are stable by this derivation.\\
Recall that we note $B[a] $ the  unitary $B$-algebra generated by $1$ and $a$ over $B$ with the commutation relation $ab - ba = b^2$. It is the free left $B$-module with basis $1, a, \dots, a^n, \dots$. Its product is defined by  the commutation relation $aS - Sa = b^2S'$ for $S \in B$. Recall also that  $\A = \widehat{B}[a]$ is  the unitary $\widehat{B}$-algebra generated by $1$ and  $a$ with the same commutation relation for $S \in \widehat{B}$.

\begin{defn}\label{(a,b) 1}
We define a {\bf convergent (a,b)-module} as a free, finite rank $B$-module $\mathcal{E}$ endowed with a continuous\footnote{for the natural topology of $B$.} $\mathbb{C}$-action of $a : \mathcal{E} \to \mathcal{E}$ such that 
$ab - ba = b^2$. Then we define its {\bf formal completion} (in $b$) of $\mathcal{E}$ as  $E := \mathcal{E}\otimes_B \widehat{B}$ on  which the action of $a$ on it is defined by
$$ a(x\otimes S) = ax\otimes S + x\otimes b^2S'$$
for $x \in \mathcal{E}$ and $S \in \widehat{B}$. It is a (formal) (a,b)-module (see \cite{[B.93]}). \\
\end{defn}

\parag{Remarks}
\begin{enumerate}
\item The continuity of $a$ implies that for any $S \in B$ we have $aS -Sa = b^2S'$ as an equality between $\mathbb{C}$-linear continuous endomorphisms of $\mathcal{E}$.\\
 Then a convergent (a,b)-module is a left $B[a]$-module.
\item A sub-module of a convergent (a,b)-module $\mathcal{E}$ is, by definition, a sub-$B$-module of $\mathcal{E}$ which is stable by $a$. So a sub-module is simply a left  sub-$B[a]$-module.\\
As any sub-$B$-module of a free finite rank $B$-module is again
free and finite rank $B$-module, a sub-module of a convergent (a,b)-module $\mathcal{E}$ is itself a convergent (a,b)-module.
\item Remark that if $\mathcal{F}$ is a sub-module of the convergent (a,b)-module $\mathcal{E}$ the quotient $\mathcal{E}/\mathcal{F}$ is not, in general, a free $B$-module, because it may have $b$-torsion. This is the reason to introduce  the notion of {\bf normal} sub-module.
\end{enumerate}

\begin{defn}\label{normal}
Let $\mathcal{F} \subset \mathcal{E}$ be a sub-module of the convergent (a,b)-module $\mathcal{E}$. We say that $\mathcal{F}$ is {\bf normal} when it satisfies  $\mathcal{F} \cap b\mathcal{E} = b\mathcal{F}$.
\end{defn}

This condition is necessary and sufficient in order that  the quotient  $\mathcal{E}/\mathcal{F}$ has no $b$-torsion and so that $\mathcal{E}/\mathcal{F}$ is a convergent (a,b)-module.

\begin{lemma}\label{normalization}
Let $\mathcal{F} \subset \mathcal{E}$ be a sub-module of the convergent (a,b)-module $\mathcal{E}$. Then define the {\bf normalization} $\tilde{\mathcal{F}}$ of $\mathcal{F}$ by the equality:
$$\tilde{\mathcal{F}} := \{ x \in \mathcal{E} / \ /  \exists n \in \mathbb{N} \quad {\rm such \ that} \quad b^nx \in \mathcal{F} \}.$$
Then $\tilde{\mathcal{F}}$ is a normal sub-module of $\mathcal{E}$ and it is the smallest normal sub-module in $\mathcal{E}$ which contains $\mathcal{F}$. The quotient $\tilde{\mathcal{F}}/\mathcal{F}$ is a finite dimensional complex vector space.
\end{lemma}

\parag{Proof} It clear that $\tilde{\mathcal{F}}$ is a $B$-sub-module (so it is free finite rank over $B$)  and its stability by the action of  $a$ is consequence the formula $b^Na = ab^N - Nb^{N+1}$ in $B[a]$. It is clearly normal. Let $\mathcal{G}$ be a normal sub-module of $\mathcal{E}$ containing $\mathcal{F}$. Now if $x$ is in $\tilde{\mathcal{F}}$ there exists an integer $n$ such that $b^nx \in \mathcal{F} \subset \mathcal{G}$. Since $\mathcal{G}$ is normal, we have $x \in \mathcal{G}$ and so $\tilde{\mathcal{F}} \subset \mathcal{G}$. Then $\tilde{\mathcal{F}}$ is the smallest normal sub-module in $\mathcal{E}$ containing $\mathcal{F}$.\\
 Let $e_1, \dots, e_k$ be a $B$-basis of $\tilde{\mathcal{F}}$. For each $j \in [1,k]$ there exists $N_j \in \mathbb{N}$ such that $b^{N_j}e_j$ is $\mathcal{F}$. Then for $N := \sup \{N_j, j \in  [1,k]\}$ we have $b^N\tilde{\mathcal{F}} \subset \mathcal{F}$. So $\tilde{\mathcal{F}}/\mathcal{F}$ is a quotient of the finite dimensional complex vector space $\tilde{\mathcal{F}}/b^N\tilde{\mathcal{F}}$, concluding the proof.$\hfill \blacksquare$\\

\begin{defn}\label{regular}
A convergent (a,b)-module $\mathcal{E}$  has a {\bf simple pole} when it satisfies $a\mathcal{E} \subset b\mathcal{E}$.\\
A convergent (a,b)-module $\mathcal{E}$  is {\bf regular} when it is a sub-module of some simple pole convergent (a,b)-module.
\end{defn}

\begin{lemma}\label{sature}
Let $\mathcal{E}$ be a regular convergent (a,b)-module. Then there exists a natural injective $B[a]$-linear map $j : \mathcal{E} \to \mathcal{E}^\sharp$ where $\mathcal{E}^\sharp$ is a simple pole convergent (a,b)-module such that  any injective $B[a]$-linear map $h : \mathcal{E} \to \mathcal{E}_1$ into a simple pole convergent (a,b)-module $\mathcal{E}_1$ factorizes (uniquely) by a $B[a]$-linear map $H : \mathcal{E}^\sharp \to \mathcal{E}_1$.\\
Moreover, the quotient $\mathcal{E}^\sharp/j(\mathcal{E})$ is a finite dimensional complex space.
\end{lemma}

\parag{Proof} Let $K := B[b^{-1}]$ and consider on  the left  $B$-module $\mathcal{E}\otimes_B K$ the $\mathbb{C}$-linear action of $a$ defined by
 $$a(x\otimes b^{-p}) := ax\otimes b^{-p} - px\otimes b^{-p+1}$$ 
 for each $p \in \mathbb{N}$. Then it is a left  $B[a]$-module and $j_0 : \mathcal{E} \to \mathcal{E}\otimes_B K$ given by $j_0(x) := x\otimes 1$ is $B[a]$-linear and injective. Define 
 $$\mathcal{E}^\sharp := \sum_{p=0}^\infty (b^{-1}a)^p j_0(\mathcal{E}) \subset \mathcal{E}\otimes_B K$$
   and define  $j$ as the induced map. First we want to prove that $\mathcal{E}^\sharp$ is a finitely generated $B$-module (obviously stable by $b^{-1}a$ so by $a$).\\
   For this purpose consider an  injective $B[a]$-linear map $h : \mathcal{E} \to \mathcal{E}_1$ into a simple pole convergent (a,b)-module $\mathcal{E}_1$. The regularity of $\mathcal{E}$ insures that there exists at least one such map. Now the simple pole assumption on $\mathcal{E}_1$ implies that $b^{-1}a$ acts on $\mathcal{E}_1$. So for each $q \in \mathbb{N}$ we may extend the $B[a]$-linear map $h$ to an injective $B[a]$-linear map  $h_q : \sum_{p=0}^q (b^{-1}a)^p j_0(\mathcal{E}) \to \mathcal{E}_1$ commuting with $b^{-1}a$. This defines an increasing sequence of $B$-sub-modules in $\mathcal{E}_1$. So it is stationary for  $q  \geq q_0$ for some integer $q_0$. Since $h_q$ is still injective for any $q$, this implies that
    $$\mathcal{E}^\sharp := \sum_{p=0}^{q_0} (b^{-1}a)^p j_0(\mathcal{E}) $$
    and $\mathcal{E}^\sharp $  is a simple pole convergent (a,b)-module containing $\mathcal{E}$ via the map $j$ induced by $j_0$.
   The previous argument shows the universal property of  the inclusion map constructed above, $j: \mathcal{E} \to \mathcal{E}^\sharp$, relative to any injective $B[a]$-linear map of $\mathcal{E}$ into  a simple pole convergent (a,b)-module.\\
 The only point to conclude the proof if to show the finite dimension of the complex vector space $\mathcal{E}^\sharp/j(\mathcal{E})$. But we already know that
   $b^{q_0}\mathcal{E}^\sharp \subset j(\mathcal{E})$ and this is enough to conclude, as $\mathcal{E}^\sharp$ is a finite $B$-module and $\mathcal{E}^\sharp/j(\mathcal{E})$ is a quotient of $\mathcal{E}^\sharp/b^{q_0}\mathcal{E}^\sharp$.$\hfill \blacksquare$
   
   \parag{Remark} Consider a  short exact sequence of regular convergent (a,b)-modules $0 \to \mathcal{F} \to \mathcal{E} \to \mathcal{G} \to 0$. It gives a surjective map $\mathcal{E}^\sharp \to \mathcal{G}^\sharp$ thanks to the minimality\footnote{The image of $\mathcal{E}^\sharp$ in $\mathcal{G}[b^{-1}]$ has a simple pole and contains $\mathcal{G}$.} of $\mathcal{G}^\sharp$. But the kernel of this map is, in general, bigger than $\mathcal{F}^\sharp$ although is has a simple pole, because a normal sub-module of a simple pole module has again a simple pole.
   
   \parag{Exercise} Show that in the situation of the previous remark the kernel of the  surjective map $\mathcal{E}^\sharp \to \mathcal{G}^\sharp$ is equal to
   $N_{\mathcal{E}^\sharp}(\mathcal{F}^\sharp) = N_{\mathcal{E}^\sharp}(\mathcal{F})$.
   
   \parag{Warning} From now on we omit ``convergent" when we consider an (a,b)-module, and we use ``script characters" (like $\mathcal{E}$) for these. If we want to consider a ``formal" (a,b)-module we shall use ``roman characters" (like $E$) and say ``formal (a,b)-module" if we want more precision.
   
   \subsection{The Bernstein polynomial of a regular (a,b)-module.}
   
   We introduce now a fundamental (numerical) invariant of a regular (a,b)-module.
   
   \begin{defn}\label{Bernstein 0}
   Let $\mathcal{E}$ be a regular (a,b)-module. The Bernstein polynomial of $\mathcal{E}$ is the minimal polynomial of the action of $-b^{-1}a$ on the finite dimensional complex  vector space $\mathcal{E}^\sharp/b\mathcal{E}^\sharp$.
   \end{defn}
   
   \parag{Remarks} \begin{enumerate}
   \item If $E$ is the formal $b$-completion of the regular (a,b)-module $\mathcal{E}$, so $E := \mathcal{E}\otimes_{B} \widehat{B}$;  then $E^\sharp$ is the formal completion of $\mathcal{E}^\sharp$ and there is a natural isomorphism $\mathcal{E}^\sharp/b\mathcal{E}^\sharp \simeq E^\sharp/bE^\sharp$ which commutes with the respective actions of $b^{-1}a$.\\
    So the Bernstein polynomial of $\mathcal{E}, \mathcal{E}^\sharp, E^\sharp$ and of $E$ are the same.
   \item Let $m$ be a non negative integer. When $\mathcal{E}$ is a simple pole (a,b)-module, the sub-module $b^m\mathcal{E}$ which has finite complex  co-dimension in $\mathcal{E}$ is again a simple pole (a,b)-module since $ab^m\mathcal{E} = b^m(a + mb)\mathcal{E} \subset  b(b^m\mathcal{E})$. The Bernstein polynomial of $b^m\mathcal{E}$ is then given by 
   $$ B_{b^m\mathcal{E}}(x) = B_{\mathcal{E}}(x+m)$$
   %since $(b^{-1}a)^q b^m = b^m(b^{-1}a + m)^q$.
     \item Let $\pi : \mathcal{E} \to \mathcal{G}$ a surjective $B[a]$-linear map. The map $\mathcal{E}^\sharp \to \mathcal{G}^\sharp$ is surjective  (see the remark at the end of section 3.1). Since this map commutes with the respective actions of $b^{-1}a$, the Bernstein polynomial of $\mathcal{G}$ divides the Bernstein polynomial of $\mathcal{E}$.
   \item Despite the previous remark, there are, in general,  two difficulties to compute the Bernstein polynomial of a regular (a,b)-module $\mathcal{E}$: \\
   The first one comes from the non left exactness of the functor ``$\ \sharp \ $" (see the remark following  Lemma \ref{sature}).\\
   The second difficulty comes from  the fact that the minimal polynomial of a pair $(V, T)$ of a vector space $V$ with an endomorphism of $V$ also does not behave nicely under injective maps compatible with $T$.\\
   This is the reason to introduce below the notion of fresco, which is stable by short exact sequences and  for which we shall dispose of a nice behavior of their Bernstein polynomials in  short exact sequences (see \cite{[B.09]}). . 
   \end{enumerate}

   We shall now prove the following important key for the Decomposition Theorem \ref{DT}.
   
   \begin{prop}\label{Sol. 0}
   Let $0 \to \mathcal{F} \to \mathcal{E} \to \mathcal{G} \to 0$ be an exact sequence of simple pole (a,b)-modules. Assume that for each root $-\lambda$ of  $B_\mathcal{F}$ and each root $-\mu$ of $B_\mathcal{G}$ we have $\mu - \lambda \not\in \mathbb{N}^*$. Then this short exact sequence splits and there exists a normal sub-module $\mathcal{G}_0$ in $\mathcal{E}$ such that $\mathcal{E} = \mathcal{F} \oplus \mathcal{G}_0$.
   \end{prop}
   
   The proof of this proposition  uses the following lemma
   
   \begin{lemma}\label{element.}
Let $F$ and $G$ two matrices with complex entries of size $(l,l)$ and $(k,k)$ respectively having no common eigenvalue. Then consider the endomorphisms $f$ and $g$ on the vector space of  matrices $Z$ with complex entries and size $(k,l)$ given by left and write multiplication by $F$ and $G$ respectively. Then $f - g$ is bijective.
\end{lemma}

\parag{proof} Since the endomorphisms $f$ and $g$ commute we may find a basis $Z_1, \dots ,Z_{kl}$ of the vector space of the  $(k,l)$ matrices which makes the matrices of $f$ and $g$ lower triangular. So for each $i \in [1,kl]$ we have $FZ_i = \lambda_i Z_i $ modulo $V_{i+1}$ and $Z_i G = \mu_iZ_i $ modulo $V_{i+1}$ where $V_{i+1}$ is the subspace generated by $Z_{i+1}, \dots, Z_{kl}$.\\
Let $Z := \sum_{i=1}^{kl}  \alpha_i Z_i$ such that $FZ = ZG$ and assume that $Z \not= 0$. Let  $i_0$ be  the smallest integer in $[1, kl]$ such that $\alpha_{i_0} \not= 0$. Then we have
$$ FZ = \alpha_{i_0}\lambda_{i_0}Z_{i_0} + V_{i_0+1} \quad {\rm and} \quad ZG = \alpha_{i_0}\mu_{i_0}Z_{i_0} + V_{i_0+1} $$
which implies that $\lambda_{i_0} = \mu_{i_0}$. But the eigenvalues of $f$ (resp. of g) are eigenvalues of $F$ (resp. of G) because $FZ = \lambda Z $ implies that each column of $Z$ is an eigenvector of $F$ (or $0$) and  if $Z$ is not zero, at least one column of $Z$ is not zero (resp. if $ZG = \mu Z$ and $Z$ is not zero, at least one line of $Z$ is a non zero eigenvector for the transpose of $G$ and $\mu$ is an eigenvalue for $G$). So we obtain a contradiction assuming that a non zero $Z$ satisfies $FZ = ZG$.$\hfill \blacksquare$\\

\parag{Proof of Proposition \ref{Sol. 0}} Let $e := (e_1, \dots, e_k)$ be a $B$-basis of $\mathcal{F}$ and let $(e_1, \dots, e_k, \varepsilon_1, \dots, \varepsilon_l)$ be a $B$-basis of $\mathcal{E}$ such we have (with matrix notations):
$$ a\begin{pmatrix} e \\ \varepsilon \end{pmatrix} = b\begin{pmatrix} F & 0 \\ bX & G\end{pmatrix}\begin{pmatrix} e \\ \varepsilon \end{pmatrix} $$
where $F, G$ and $X$ are respectively $(k,k), (l,l)$ and $(l,k)$ matrices with entries in $B$. This possible because we have a direct sum decomposition of $\mathcal{E}/b\mathcal{E}$ compatible with the spectral decomposition of the action of $b^{-1}a$ on this finite dimensional  vector space. Then  write: 
\begin{equation*}
 F := \sum_{j=0}^\infty  F_jb^j,  \  G := \sum_{j=0}^\infty  G_jb^j,  \  X := \sum_{j=0}^\infty  X_jb^j 
\end{equation*}
where $F_j, G_j, X_j$ are matrices with complex coefficients. Choose on the vector spaces of matrices with complex entries and size $(k,k), (l,l), (l,k)$ norms such that when the product is defined we have $\vert XY\vert \leq \vert X\vert\times \vert Y\vert$. Choose also a norm on the vector space of endomorphism of $(l,k)$-matrices such that $\vert H(Z)\vert \leq \vert\vert H\vert\vert \times \vert Z\vert$.\\
Note that the endomorphism $H$ defined by $H(Z) = ZF_0 - G_0Z$  has no eigenvalue in $-\mathbb{N}^*$ thanks to our hypothesis and the lemma above applied to $H + j Id, j \in \mathbb{N}^*$.\\
Now we look for a $(l,k)$ matrix $Z$ with entries in $B$ such that we have
$$a(\varepsilon + Ze) = bG(\varepsilon + Ze).$$ 
Put $Z := \sum_{j=1}^\infty Z_jb^j $. Then to find $Z$ with entries in $\widehat{B}$ is equivalent to solve the equation
$$ a(\varepsilon + Ze) = b^2Xe + bG\varepsilon + bZFe + b^2Z'e = bG\varepsilon + bGZe $$
which is equivalent to 
$$ ZF -GZ +bZ' = -bX$$
and then  is equivalent to  the system of equations
\begin{equation*}
Z_jF_0 - G_0Z_j + jZ_j = -X_{j-1} + \sum_{p=1}^j \big(G_pZ_{j-p} - Z_{j-p}F_p\big) \quad \forall j \in \mathbb{N}^*. \tag{S}
\end{equation*}
Since the endomorphism $H + j Id$ is bijective for each $j \geq 1$, an induction shows that there exists a unique solution $Z$ with entries in $\widehat{B}$.\\
We want now to show that $Z$ has its entries in $B$. So fix $R > 1$ and choose a positive constant $C_R$ such that we have the estimates:
\begin{equation}
 \vert F_j\vert \leq C_RR^j j! \quad \vert G_j\vert \leq C_RR^j j! \quad \vert X_j\vert \leq C_RR^j j! 
\end{equation}
Let $\rho := \vert\vert H\vert\vert$. Then for $j > \rho$ we have
\begin{equation}
 (H + j Id)^{-1} = j^{-1}(Id + H/j)^{-1} =  j^{-1} \sum_{h=0}^\infty  (-H/j)^h 
 \end{equation} and this gives the estimates
 \begin{equation}
 \vert\vert (H + j.Id)^{-1} \vert\vert \leq \frac{1}{j - \rho}  \quad \forall j > \rho.
 \end{equation}
Then choose a positive constant $D_R \geq 1$ large enough such that the estimate  \\
$\vert Z_j \vert \leq D_RR^j j!$  \ is valid  for any $j \leq \rho + 5C_R$.\\
Now assume that for some $j_0 \geq \rho + 5C_R$ we have proved that $\vert Z_j \vert \leq D_RR^j j!$ for any $j \leq j_0 -1$. We shall prove that this estimate is also valid for $j_0$.\\
We have
$$  H(Z_{j_0}) + j_0Z_{j_0} =  -X_{j_0-1} + \sum_{p=1}^{j_0} G_pZ_{j_0-p} - Z_{j_0-p}F_p  $$
and so, using the estimate $(9)$ for the norm of  $(H + j Id)^{-1}$:
$$ \vert Z_{j_0} \vert \leq \frac{1}{j_0 -\rho}\big(C_R R^{j_0-1}(j_0-1)! + 2C_RD_R R^{j_0} (j_0)! \sum_{p=1}^{j_0} \frac{(j_0 -p)! p!}{j_0!}\big) $$
We obtain, using $j_0 \geq \rho + 5C_R$, \  $D_R \geq 1/R$ and Lemma \ref{inegalite 1}:
$$  \vert Z_{j_0} \vert \leq \frac{1}{j_0 -\rho} C_R R^{j_0} j_0!\big( \frac{1}{R j_0} + 4 D_R\big) \leq D_R R^{j_0} j_0! $$
since 
$$ \frac{1}{j_0 -\rho}\big( \frac{1}{R j_0} + 4 D_R\big) \leq 5D_R/5C_R.$$
This completes the proof that $Z$ has its entries in $B$.\\
The conclusion follows immediately defining $\mathcal{G}_0$ as the sub-module generated by $\varepsilon + Ze$.$\hfill \blacksquare$

\parag{Remark} The uniqueness of the matrix $Z$ in the proposition above implies that the splitting of the exact sequence is unique in this situation. This means that the complement $\mathcal{G}_0$ constructed in the proof is unique and that the decomposition $\mathcal{E} \simeq \mathcal{F} \oplus \mathcal{G}_0$ is ``natural".\\

\begin{cor}\label{Sol. 1}
 Let $\mathcal{E}$ be a simple pole convergent  $(a,b)$--module with basis  the column $e := \, ^t(e_1, \dots, e_k)$ with $ae = bF(b)e$ where $F$ is a $(k,k)$-matrix with entries in  the algebra 
    $B = \mathbb{C}\{\{b\}\}$.  Write $F(b) := \sum_{j=0}^{+\infty} F_jb^j$ and assume that, for a given complex number $\lambda $, the spectrum of $F_0$ does not meet $\{ \lambda - \mathbb{N}\}$. Then for any $y \in \mathcal{E}$ there exists a unique $x \in \mathcal{E}$ such that
\begin{equation} 
(a - \lambda b)x = by \quad {\rm in} \ \mathcal{E} 
\end{equation}
\end{cor}

\parag{Proof} Let $e$ be a $B$-basis of $\mathcal{E}$ which satisfies $ae = bFe$ with $F := \sum_{j=0}^\infty  F_jb^j$ is a $(k,k)$-matrix with entries in $B$. Then for a given $= Ye, Y \in B^k$, we to look for  $Z \in B^k$ such that  $x := Ze$, satisfies  $(a-\lambda b)x = by$. This leads to the equation 
$$ ZbFe + b^2Z'e - \lambda bZe = bYe $$
which is equivalent to the system of equations, writng $Z = \sum_{j=0}^\infty Z_jb_j$:
$$ Z_jF_0 - \lambda Z_j + jZ_j = Y_{j} - \sum_{p=1}^j Z_{j-p}F_p \quad \forall j \geq 0 $$
For $j \geq 1$ this is the same system as the system $(S)$  in the previous proof  by letting $G = G_0 = \lambda Id$. As the equation for $j=0$ has a unique solution, since our hypothesis implies that $F_0 -\lambda $ is bijective  on $\mathcal{E}/b\mathcal{E}$, we find a unique solution $Z$ with entries in $B$ for each $y \in \mathcal{E}$ thanks to the estimates given in the proof of Proposition \ref{Sol. 0}.$\hfill \blacksquare$

\parag{Remark} The case $\lambda = 0$ of the previous corollary shows that a simple pole (a,b)-module $\mathcal{E}$ such that its Bernstein polynomial has no root in 
$\mathbb{N}$ satisfies $a\mathcal{E} = b\mathcal{E}$, so $b^{-1}a$ is bijective. \\

\begin{cor}\label{sol. Jordan}
Let $\mathcal{E}$ be a simple pole convergent (a,b)-module and let $\lambda$ be an eigenvalue of $b^{-1}a$ acting on $\mathcal{E}/b\mathcal{E}$ and having the following properties:
\begin{enumerate}
\item $\lambda -p$ is not an eigenvalue of $b^{-1}a$ acting on $\mathcal{E}/b\mathcal{E}$ for each $p \in \mathbb{N}^*$.
\item There exists a rank $k$ Jordan bloc for the eigenvalue $\lambda$ acting on $\mathcal{E}/b\mathcal{E}$
\end{enumerate}
Then there exists elements $\varepsilon_1, \dots, \varepsilon_k$ in $\mathcal{E}$ which are independent over $B$ and satisfies the relations
\begin{equation}
a\varepsilon_j = \lambda b\varepsilon_j + b\varepsilon_{j+1} \quad {\rm for \ each} \ j \in [1, k] \quad {\rm with \ the \ convention} \quad \varepsilon_{k+1} = 0.
\end{equation}
\end{cor}

\parag{Proof} Thanks to Corollary \ref{Sol. 1} we can solve, for any $y \in \mathcal{E}$, the equation $(a - (\lambda-1)b)x = by$. \\
Then consider elements $e_1, \dots, e_k$ in $\mathcal{E}$ such that they induce a $k$-Jordan bloc for the eigenvalue $\lambda$ of $b^{-1}a$ acting on $\mathcal{E}/b\mathcal{E}$.
Then there exist $y_1, \dots, y_k$ in $\mathcal{E}$ such the following relations holds
\begin{align*}
& ae_1 = \lambda be_1 +  be_2 + b^2y_1 \\
&ae_2 = \lambda be_2 + be_3 + b^2y_2 \\
& \dots \dots \\
& ae_k = \lambda be_k + b^2y_k
\end{align*}
We look now for $x_1, \dots, x_k$ in $\mathcal{E}$ such that $\varepsilon_j := e_j - bx_j$ satisfy the equations $(11)$. We shall argue by a descending induction on $j$. Assume that we have already found $x_{j+1}, \dots, x_{k+1} := 0$ for some $j \in [1, k]$. Then to find $x_j$ we have to solve the equation
$$ (a - \lambda b)(e_j - bx_j) = b(e_{j+1} - bx_{j+1}) $$
which is equivalent, thanks to the relation $(a-\lambda b)b = b(a - (\lambda-1)b$ and the injectivity of $b$,  to
$$ (a - (\lambda-1)b)x_j = b(y_j - x_{j+1}).$$
This is enough to conclude the proof.$\hfill \blacksquare$

\parag{Remark} An easy consequence of the previous result is that if $-\alpha$ is the biggest root of $B_\mathcal{E}$ for a regular (a,b)-module $\mathcal{E}$ there exists $x \in \mathcal{E}^\sharp \setminus b\mathcal{E}^\sharp $ such that $ax = \alpha bx$. Compare with Lemma \ref{Bernstein 2} below.$\hfill \square$\\

 The following  corollary is immediate.
 
\begin{cor}\label{class. 1}
A rank 1 regular convergent (a,b)-module is isomorphic to a quotient $B[a]/B[a].(a-\lambda b) $. So it has a $B$-basis $e_\lambda $, and  $\mathcal{E}_\lambda := Be_\lambda$ where $ae_\lambda = \lambda be_\lambda$.\\ 
It has a simple pole, the action of  $b^{-1}a$ on the  $1$-dimensional vector space $\mathcal{E}/b\mathcal{E}$ is given by $\lambda$ and its Bernstein polynomial is $x + \lambda$.$\hfill \blacksquare$
\end{cor}

Note  that the previous proposition gives also that if $\mathcal{E}$ is a simple pole rank $k$ (a,b)-module such that $\mathcal{E}/b\mathcal{E}$ has an unique Jordan block which has  rank $k$ and  eigenvalue $\lambda$, then $\mathcal{E}$ is isomorphic to the (a,b)-module with $B$-basis $\varepsilon_1, \dots, \varepsilon_k$ where the action of $a$ is defined by the relations $(11)$ above. \\

The following application of  Corollary   \ref{sol. Jordan} is rather useful.

\begin{lemma}\label{big 0}
Let $\mathcal{F}$ be a submodule of a regular (a,b)-module $\mathcal{E}$. If $-\beta$ is a root of the Bernstein polynomial of $\mathcal{F}$ there exists an integer $m \in \mathbb{N}$ such that $-\beta + m$ is a root of the Bernstein polynomial of $\mathcal{E}$.
\end{lemma}

\parag{Proof}  It is enough to prove the result for the inclusion of $\mathcal{F}^\sharp$ in $\mathcal{E}^\sharp$. Moreover, we may assume that $-\beta$ is the biggest root of the Bernstein polynomial of $\mathcal{F}^\sharp$ in $-\beta + \mathbb{Z}$. So $\beta$ is an eigenvalue of the action of $b^{-1}a$ on $\mathcal{F}^\sharp\big/b\mathcal{F}^\sharp$ which is minimal in $\beta  + \mathbb{Z}$. Then  Corollary \ref{sol. Jordan} gives the existence of an $x$ in $\mathcal{F}^\sharp \setminus b\mathcal{F}^\sharp$ which satisfies $(a - \beta b)x = 0$. Let $p \in \mathbb{N}$ be the maximal integer such that $x$ is in $b^p\mathcal{E}^\sharp$ and write $x = b^py$ with $y \in \mathcal{E}^\sharp$. Then $y$ satisfies $(a - (\beta -p)b)y = 0$ and $y \not\in b\mathcal{E}^\sharp$ by definition of $p$. So we see that $-\beta + p$ is a root of the Bernstein polynomial of $\mathcal{E}$ concluding the proof.$\hfill \blacksquare$\\

\subsection{The Decomposition Theorem}

We begin by an easy variant of  Lemma \ref{big 0} which  is a useful tool to discuss $\mathscr{A}$-primitive (a,b)-modules.

\begin{lemma}\label{Bernstein 2}
Let $\mathcal{E}$ be a regular (a,b)-module and assume that $-\lambda$ is a root of its Bernstein polynomial. Then there exists $m \in \mathbb{N}$ and a non zero $x$ in $\mathcal{E}$ such that
$(a - (\lambda + m)b)x = 0$. Conversely, if there exists a non zero $x \in \mathbb{E}$ such that $(a - \lambda b)x = 0$ then  $B_\mathcal{E}$ has a root in $-\lambda + \mathbb{N}$.
\end{lemma}

\parag{Proof} Let $-\lambda_1 := -\lambda +m_1$ be the biggest root of $B_{\mathcal{E}}$ which is in $-\lambda + \mathbb{N}$ and fix  $y \in \mathcal{E}^\sharp \setminus b\mathcal{E}^\sharp$ which satisfies
$(a - \lambda_1b)y \in b^2\mathcal{E}^\sharp$. Then put $(a - \lambda_1b)y = b^2z$. Since there is no root of $B_{\mathcal{E}}$ in $-\lambda_1+1+\mathbb{N}$, we know (see Corollary \ref{sol. Jordan}) that $b^{-1}a - (\lambda_1-1)$ is bijective on $\mathcal{E}^\sharp$. Then there exists $x \in \mathcal{E}^\sharp$ satisfying  $(b^{-1}a -(\lambda_1-1))x = z$ and so  we have $b(a - (\lambda_1-1) b)x = (a -\lambda_1 b) bx = b^2z = (a - \lambda_1 b)y $ and  $(a - \lambda_1 b)(y - bx) = 0$. Moreover, as $y$ is not in $b\mathcal{E}^\sharp$, $y -bx$ is also not in $b\mathcal{E}^\sharp$ and then $y - bx \not= 0$.\\
Let $n \in \mathbb{N}$ such that $b^n\mathcal{E}^\sharp \subset \mathcal{E}$ and $n \geq m_1$. Then $t := b^n(y - bx)$ is in $\mathcal{E} \setminus \{0\}$ and satisfies 
$(a - (\lambda - m_1 + n)b)t = 0$ with $n -m_1 \in \mathbb{N}$.\\
Conversely if $x \in \mathcal{E} \setminus \{0\}$ satisfies $(a - \lambda b)x = 0$, write $x = b^ny$ for some $y \in \mathcal{E}^\sharp \setminus b\mathcal{E}^\sharp$ (with  some $n \in \mathbb{N}$). Then $-\lambda + n$ is a root of $B-\mathcal{E}$.$\hfill \blacksquare$\\

\parag{Remark}It is easy to deduce from the previous lemma that any non zero regular (a,b)-module admits a rank $1$ normal (a,b)-sub-module. Then an easy induction gives the existence, for any regular (a,b)-module $\mathcal{E}$,  of a Jordan-H\"{o}lder sequence
 $$0 \to \mathcal{F}_1 \to \mathcal{F}_2 \to \dots \to \mathcal{F}_k = \mathcal{E}$$
 where, for each $j \in [1, k]$ with $k$ the rank of $\mathcal{E}$,  $\mathcal{F}_j$ is normal in $\mathcal{E}$ and such that $\mathcal{F}_j\big/\mathcal{F}_{j-1}$ has rank $1$ (with $\mathcal{F}_0 = \{0\}$). See \cite{[B.09]} in the formal case.\\

 Let $\mathscr{A}$ a subset of $ \mathbb{C}/\mathbb{Z} $. 
   
   \begin{defn}\label{primitive 0}
   We say that a regular (a,b)-module $\mathcal{E}$ is {\bf $\mathscr{A}$-primitive} when all roots of its Bernstein polynomial $B_{\mathcal{E}}$ are in $-\mathscr{A}+ \mathbb{Z}$.\\
   In the case where $\mathscr{A} = \{\alpha\}$ where $\alpha$ is an element of $\mathbb{C}/\mathbb{Z}$  we say that $\mathscr{E}$ is $[\alpha]$-primitive.   
   \end{defn}

The next proposition is a first step to the Decomposition Theorem \ref{DT}.

  \begin{prop}\label{primitive 4}
   Let $\mathscr{A}$ be a subset in $\mathbb{C}/\mathbb{Z}$ and let $\mathcal{E}$ be a regular (a,b)-module. Then there exists a maximal sub-module $\mathcal{E}_{[\mathscr{A}]}$ in $\mathcal{E}$ which is $\mathscr{A}$-primitive. Moreover, this sub-module is normal and the quotient  $\mathcal{E}/\mathcal{E}_{[\mathscr{A}]}$ is a $\mathscr{A}^c$-primitive, where $\mathscr{A}^c$ is the complement of $\mathscr{A}$ in $\mathbb{C}/\mathbb{Z}$.
   \end{prop}

The proof of this result needs some lemmas.

\begin{lemma}\label{primitive 2}
Fix any subset $\mathscr{A}$ in $\mathbb{C}/\mathbb{Z}$. Let $0 \to \mathcal{F} \to \mathcal{E} \to \mathcal{G} \to 0 $ be a short exact sequence of regular (a,b)-modules. Then $\mathcal{E}$ is $\mathscr{A}$-primitive if and only if both $\mathcal{F}$ and $\mathcal{G}$ are $\mathscr{A}$-primitive.
\end{lemma}

\parag{Proof} First assume that $\mathcal{E}$ is $\mathscr{A}$-primitive. Let $-\mu$ a root of $B_{\mathcal{F}}$. Then thanks to Lemma \ref{Bernstein 2}  there exists an integer $m$ and a non zero element $x \in \mathcal{F}$ satisfying $(a - (\mu + m)b)x = 0$. Let $k$ the maximal integer such that $x$ is in $b^k\mathcal{E}^\sharp$ and put $x = b^k y$ where $y \in \mathcal{E}^\sharp$. Then we have
$(a - (\mu+m-k)b)y = 0$ with $y \not\in b\mathcal{E}^\sharp$ so $[\mu]$ is in $\mathscr{A}$. Then $\mathcal{F}$ is $\mathscr{A}$-primitive.\\
The fact that $\mathcal{G}$ is also $\mathscr{A}$-primitive is a consequence of  Remark 3 following Definition \ref{Bernstein 0}.\\
Let assume now that $\mathcal{F}$ and $\mathcal{G}$ are $\mathscr{A}$-primitive and let $-\mu$ be a root of $B_{\mathcal{E}}$. Assume that $[\mu] \not\in \mathscr{A}$. There exists a non zero element $x \in \mathcal{E}$ which satisfies $(a - (\mu-m)b)x = 0$, where $m$ is an integer, thanks to Lemma \ref{Bernstein 2}. Then the same argument  implies that the image of $x$ in $\mathcal{G}$ is zero. So $x$ is in $\mathcal{F}$ which contradicts  the non vanishing of $x$.$\hfill \blacksquare$

\begin{lemma}\label{primitive 3}
Let $\mathcal{E} \subset \mathcal{E}'$ be an inclusion of two regular (a,b)-modules such that $\mathcal{E}'/\mathcal{E}$ is a finite dimensional complex vector space. Then for any subset $\mathscr{A}$ in $\mathbb{C}/\mathbb{Z}$,  $\mathcal{E}'$ is $\mathscr{A}$-primitive if and only if $\mathcal{E}$ is $\mathscr{A}$-primitive.\\
\end{lemma}

\parag{Proof} Thanks to Lemma \ref{Bernstein 2} and to Corollary \ref{primitive 1} it is enough to prove the equivalence of  the existence in $\mathcal{E}$ or in $\mathcal{E}'$   for any $[\mu] \in \mathbb{C}/\mathbb{Z}$  of a  non zero solution of the equation $(a - (\mu+m)b)x = 0$ for each $m \in \mathbb{Z}$. But this is an obvious consequence of the existence of an integer $N$ such that $b^N\mathcal{E}' \subset \mathcal{E}$.$\hfill \blacksquare$\\

 \parag{Proof of Proposition \ref{primitive 4}}  Then let us prove that  the sum of two $\mathscr{A}$-primitive sub-modules  is again  $\mathscr{A}$-primitive. This is obvious for a direct sum, and in general, the sum is a quotient of the direct sum, so we conclude using Lemma \ref{Bernstein 2}. Let  $\mathcal{F}$ be a maximal $\mathscr{A}$-primitive sub-module of $\mathcal{E}$. This exists thanks to the Noether's property of $B$. Since $\mathcal{F}$ has finite co-dimension in its normalization, Lemma \ref{primitive 3} implies that $\mathcal{F}$ is normal.\\
   Assume that the Bernstein polynomial  of $\mathcal{E}/\mathcal{F}$ has a root $-\mu$ in $-\mathscr{A}$, then Lemma \ref{Bernstein 2} gives us a non zero $z \in \mathcal{E}/\mathcal{F}$ such that $(a - (\mu+m)bz = 0$ for an integer $m \in \mathbb{Z}$. Let $\mathcal{G} := B z \subset \mathcal{E}/\mathcal{F}$. It is a rank $1$  $\mathscr{A}$-primitive sub-module and its pull-back in $\mathcal{E}$ is 
   a $\mathscr{A}$-primitive sub-module, thanks to Lemma \ref{primitive 2}, which is strictly bigger than $\mathcal{F}$. Contradiction. So $\mathcal{E}/\mathcal{F}$ is $\mathscr{A}^c$-primitive.$\hfill \blacksquare$\\

\begin{cor}\label{5/2-a}
Let $\mathcal{E} \subset \mathcal{E}'$ be an inclusion of two regular (a,b)-modules. Then we have $\mathcal{E}_{[\mathcal{A}]} = \mathcal{E} \cap \mathcal{E}'_{[\mathcal{A}]} $.
\end{cor} 

\parag{Proof} We first consider the case where $\mathcal{E}'/\mathcal{E}$ is a finite dimensional complex vector space. Then $\mathcal{E} \cap \mathcal{E}'_{[\mathcal{A}]}$ has finite codimension in 
$ \mathcal{E}'_{[\mathcal{A}]}$ which is $\mathcal{A}$-primitive, so, thanks to   Lemma \ref{primitive 3}, it is $\mathcal{A}$-primitive and then contained in 
$\mathcal{E}_{[\mathcal{A}]}$. \\
Now, since $b^n\mathcal{E}' \subset \mathcal{E}$,  the previous reasoning gives $\mathcal{E}_{[\mathcal{A}]} \cap b^n \mathcal{E}' \subset   (b^n \mathcal{E}')_{[\mathcal{A}]}$.  The equality 
$(b\mathcal{F})_{[\mathcal{A}]} = b(\mathcal{F}_{[\mathcal{A}]})$ for any $\mathcal{F}$ gives then
$$ b^n \mathcal{E}_{[\mathcal{A}]} \subset \mathcal{E}_{[\mathcal{A}]} \cap b^n \mathcal{E}' \subset   (b^n \mathcal{E}')_{[\mathcal{A}]} = b^n\mathcal{E}'_{[\mathcal{A}]} $$
concluding the proof of this first case.\\
Consider now  the case where $\mathcal{F}$ is normal in $\mathcal{E}$. Then  $\mathcal{F} \cap \mathcal{E}_{[\mathcal{A}]}$ is normal in $ \mathcal{E}_{[\mathcal{A}]}$ and so its Bernstein polynomial has only roots in $-\mathcal{A}$ which gives the inclusion
$$\mathcal{F} \cap \mathcal{E}_{[\mathcal{A}]} \subset  \mathcal{F}_{[\mathcal{A}]} .$$
Now if $x$ is in $\mathcal{F}_{[\mathcal{A}]}$ then $N_{\mathcal{E}}(B[a]x)$ is contained in $\mathcal{F} \cap \mathcal{E}_{[\mathcal{A}]}$ since $B[a]x$ has finite codimension in $N_{\mathcal{E}}(B[a]x)$  and is $\mathcal{A}$-primitive\footnote{As a sub-module of the $\mathcal{A}$-primitive module $\mathcal{F}_{[\mathcal{A}]}$.}, so we may apply  the first part of the proof  and  we have the equality
$$\mathcal{F} \cap \mathcal{E}_{[\mathcal{A}]} =  \mathcal{F}_{[\mathcal{A}]} $$
when $\mathcal{F}$ is normal. \\
When $\mathcal{F}$ is not assumed to be normal,  $\mathcal{F}$  has finite co-dimension in its normalisation $\mathcal{G} := N_{\mathcal{E}}(\mathcal{F})$ and we have
$$ \mathcal{F}_{[\mathcal{A}]} =  \mathcal{F} \cap   \mathcal{G}_{[\mathcal{A}]} = \mathcal{F} \cap\big(\mathcal{G}  \cap \mathcal{E}_{[\mathcal{A}]}\big) =  \mathcal{F} \cap \mathcal{E}_{[\mathcal{A}]}$$
by the first two cases, and the conclusion follows.$\hfill\blacksquare$

\parag{Remark} The previous corollary applies to the inclusion $\mathcal{E} \subset \mathcal{E}^\sharp$ and gives for any $\mathcal{A} \subset \mathbb{C}\big/\mathbb{Z}$ the equality $(\mathcal{E}^\sharp)_{[\mathcal{A}]} \cap \mathcal{E} = \mathcal{E}_{[\mathcal{A}]}$. So $\mathcal{E}_{[\mathcal{A}]}$ has finite co-dimension in $(\mathcal{E}^\sharp)_{[\mathcal{A}]}$ and $\mathcal{E}_{[\mathcal{A}]}, (\mathcal{E}_{[\mathcal{A}]})^\sharp $ and $(\mathcal{E}^\sharp)_{[\mathcal{A}]}$ have the same rank.\\

\begin{defn}\label{quotient 1}
For $\mathcal{E}$ a regular (a,b)-module and for $\mathcal{A} \subset \mathbb{C}\big/\mathbb{Z}$ define
$$ \mathcal{E}^{[\mathcal{A}]} := \mathcal{E}\big/\mathcal{E}_{[\mathcal{A}^c]} $$
where $\mathcal{A}^c$ is the complement of $\mathcal{A}$ in $\mathbb{C}\big/\mathbb{Z}$.
\end{defn}

Note that $\mathcal{E}_{[\mathcal{A}^c]} $ is normal in $\mathcal{E}$ thanks to Proposition \ref{primitive 4}, so that $ \mathcal{E}^{[\mathcal{A}]} $ is again a regular (a,b)-module. Note also that for a simple pole (a,b)-module $\mathcal{E}$ we have $\mathcal{E}_{[\mathcal{A}]} \simeq \mathcal{E}^{[\mathcal{A}]}$.

\begin{lemma}\label{quotient 2}
Keep the notations above and let $\mathcal{F}$ be a normal sub-module of $\mathcal{E}$ such that $\mathcal{E}\big/\mathcal{F}$ is $\mathcal{A}$-primitive. Then $\mathcal{F}$ contains $\mathcal{E}_{[\mathcal{A}^c]}$. So $\mathcal{E}\big/\mathcal{F}$ is a quotient of $\mathcal{E}^{[\mathcal{A}]}$.
\end{lemma}

\parag{proof} Let $\mathcal{G}$ be the normalization of $\mathcal{F}$ in $\mathcal{E}^\sharp$. The exact sequence
$$ 0 \to \mathcal{G} \to \mathcal{E}^\sharp \to (\mathcal{E}\big/\mathcal{F})^\sharp \to 0 $$
shows that, since $\mathcal{E}\big/\mathcal{F}$ is $[\mathcal{A}]$-primitive, a root in $-\mathcal{A}^c$ of $B_\mathcal{E}$ is a root of $B_\mathcal{G}$. Moreover
$(\mathcal{E}^\sharp)_{[\mathcal{A}^c]} \subset \mathcal{G}$ so $(\mathcal{E}^\sharp)_{[\mathcal{A}^c]}\cap \mathcal{E} \subset \mathcal{G}\cap \mathcal{E} = \mathcal{F}$ as $\mathcal{F}$ is normal in $\mathcal{E}$. So $(\mathcal{E}^\sharp)_{[\mathcal{A}^c]} \subset \mathcal{F}$ since $(\mathcal{E}^\sharp)_{[\mathcal{A}^c]} \cap \mathcal{E} = \mathcal{E}_{[\mathcal{A}^c]} $ thanks to the remark above. the conclusion follows.$\hfill\blacksquare$\\

This lemma shows that $\mathcal{E}^{[\mathcal{A}]}$ is the maximal $[\mathcal{A}]$-primitive quotient of $\mathcal{E}$. We shall call $\mathcal{E}^{[\mathcal{A}]}$  {\bf the $[\mathcal{A}]$-primitive quotient} of $\mathcal{E}$.

\begin{lemma}\label{24/2}
With the  previous notations the natural $B[a]$-linear map $\mathcal{E}_{[\mathcal{A}]} \to \mathcal{E}^{[\mathcal{A}]}$ is injective and has finite $\mathbb{C}$-dimensional co-kernel.
\end{lemma}

\parag{Proof} This map is clearly injective since $ \mathcal{E}_{[\mathcal{A}]} \cap \mathcal{E}_{[\mathcal{A}^c]} = \{0\}$. The quotient of $\mathcal{E}$ by the sub-module
$\mathcal{E}_{[\mathcal{A}]} \oplus \mathcal{E}_{[\mathcal{A}^c]} $ has finite $\mathbb{C}$-dimension since $\mathcal{E}^\sharp$ is the direct sum of
$(\mathcal{E}^\sharp)_{[\mathcal{A}]}  \oplus  (\mathcal{E}^\sharp)_{[\mathcal{A}^c]} $ and since $(\mathcal{E}_{[\mathcal{A}]})^\sharp$
 (resp. $ (\mathcal{E}_{[\mathcal{A}^c]})^\sharp $) has finite co-dimension in $(\mathcal{E}^\sharp)_{[\mathcal{A}]} $ (resp. in $  (\mathcal{E}^\sharp)_{[\mathcal{A}^c]} $).$\hfill \blacksquare$\\

 \begin{cor}\label{4/2-b}
For any regular (a,b)-module we have
$$ (\mathcal{E}^\sharp)^{[\mathcal{A}]} = (\mathcal{E}^{[\mathcal{A}]})^\sharp $$
\end{cor}

\parag{Proof}Since  $\mathcal{E}_{[\mathcal{A}^c]} $ has finite co-dimension in $(\mathcal{E}^\sharp)_{[\mathcal{A}^c]}$ and since  they have the same rank, we have 
$$N_{\mathcal{E}^\sharp}(\mathcal{E}_{[\mathcal{A}^c]}) = N_{\mathcal{E}^\sharp}((\mathcal{E})^\sharp)_{[\mathcal{A}^c]}) = (\mathcal{E}^\sharp)_{[\mathcal{A}^c]}  $$
and we obtain the desired  equality by comparison of the exact sequences
\begin{align*}
& 0 \to \mathcal{E}_{[\mathcal{A}^c]} \to \mathcal{E} \to \mathcal{E}^{[\mathcal{A}]} \to 0 \\
& 0 \to N_{\mathcal{E}^\sharp}(\mathcal{E}_{[\mathcal{A}^c]}) \to  \mathcal{E}^\sharp \to  (\mathcal{E}^{[\mathcal{A}]})^\sharp \to 0  \\
& 0 \to (\mathcal{E}^\sharp)_{[\mathcal{A}^c]} \to \mathcal{E}^\sharp \to (\mathcal{E}^\sharp)^{[\mathcal{A}]} \to 0 \qquad \qquad \qquad \qquad  \qquad \qquad  \blacksquare
\end{align*}

\bigskip

  \begin{thm}\label{DT}
  Let $\mathcal{E}$ be a simple pole  (convergent) (a,b)-module and note $-\mathscr{A}$ the image in $\mathbb{C}/\mathbb{Z}$ of the roots of the Bernstein polynomial of $\mathcal{E}$. Then 
  we have a natural isomorphism of (a,b)-modules(so of left $B[a]$-modules):
  $$ \mathcal{E} \simeq  \oplus_{\alpha \in \mathscr{A}}  \ \mathcal{E}_{[\alpha]}.$$
  \end{thm}
  
  The proof is an obvious  consequence of the following  Corollary \ref{primitive 5}  of Proposition \ref{Sol. 0}, using an induction on the cardinal of the set $\mathscr{A}$.$\hfill \blacksquare$\\
  
  For the ``naturality" of this decomposition, see the remark following Proposition  \ref{Sol. 0}.
  
  \begin{cor}\label{primitive 5}
  Let \  $0 \to \mathcal{F} \to \mathcal{E} \to \mathcal{G} \to 0 $ \  be a short exact sequence of regular (a,b)-modules. Assume that $\mathcal{E}$ has a simple pole and that $\mathcal{F}$ and
   $\mathcal{G}$ are respectively $\mathscr{A}$ and $\mathscr{B}$-primitive with $\mathscr{A} \cap \mathscr{B} = \emptyset$ in $\mathbb{C}/\mathbb{Z}$. Then the exact sequence splits and $\mathcal{E} \simeq \mathcal{F} \oplus \mathcal{G}$.$\hfill \blacksquare$\\
  \end{cor}
  
  Note that for a regular (a,b)-module  $\mathcal{E}$ which is not a simple pole (a,b)-module, such a short exact sequence does not split in general. For instance, consider  the rank $2$ (a,b)-module $\mathcal{H}$  with $B$-basis $e_1, e_2$ on which the action of $a$ is defined by
  $$ ae_1 = \lambda be_1 + e_2 \quad  ae_2 = \mu be_2 \quad  {\rm with} \quad \lambda - \mu \not\in \mathbb{Z} .$$
  It is easy to see that $\mathcal{H}_{[\lambda]} \simeq \mathcal{E}_{\lambda+1}$ is the normal sub-module  generated by the element $e_2 + (\lambda - \mu  +1)be_1$ in $\mathcal{H}$ and that $\mathcal{H}_{[\mu]}$ is generated by $e_2$. But $e_1$ is not isomorphic to $\mathcal{H}_{[\lambda]} \oplus \mathcal{H}_{[\mu]}$ since he does not have a simple pole.\\

The decomposition Theorem implies the following obvious  decomposition result for any  simple pole  (a,b)-module.

\begin{cor}\label{DT Bernstein 1}
Let $\mathcal{E}$ a simple pole (a,b)-module and let $-\mathscr{A}$ be the image in $\mathbb{Q}\big/\mathbb{Z}$ of the set of roots of its Bernstein polynomial of $\mathcal{E}$.
Then we have
$$ B_\mathcal{E} = \prod_{\alpha \in \mathscr{A}} B_{\mathcal{E}_{[\alpha]}}. \qquad \qquad \qquad \qquad  \qquad \qquad  \qquad \qquad   $$
\end{cor}

\parag{Remark} Since, for any regular (a,b)-module the map $(\mathcal{E}^\sharp)\to  (\mathcal{E}^{[\mathcal{A}]})^\sharp $ is always surjective,  the Bernstein polynomial of 
$\mathcal{E}^{[\mathcal{A}]}$ always divides $B_{\mathcal{E}}$. Moreover, the corresponding quotient has no root in $-\mathcal{A}$ so $B_{\mathcal{E}^{[\mathcal{A}]}}$ is obtained from $B_{\mathcal{E}}$ by deleting all roots which are not in $-\mathcal{A}$.\\
Note  that, in general, $B_{\mathcal{E}_{[\mathcal{A}]}}$ does not divide $B_\mathcal{E}$  as it is shown by the example below.

\parag{An example} Consider the sub-module $\mathcal{E}$  generated by $\varepsilon := e_\alpha + e_\beta$ in $E_\alpha \oplus E_\beta$. Then we have
$$ \mathcal{E}_{[\alpha]} = E_{\alpha+1} \quad {\rm and} \quad \mathcal{E}^{[\alpha]} \simeq E_\alpha $$
because $\mathcal{E}_{[\beta]} = E_{\beta+1}$ which is normal in $\mathcal{E}$, and, because  in the quotient $\mathcal{E}^{[\alpha]} := \mathcal{E}\big/\mathcal{E}_{[\beta]} $
we have $a\varepsilon = ae_\alpha = \alpha be_\alpha$ and $b\varepsilon = be_\alpha$. So the Bernstein polynomial of $\mathcal{E}^{[\alpha]}$ divides $B_{\mathcal{E}}$ and it is not the case for $\mathcal{E}_{[\alpha]}$ (which is not normal in $\mathcal{E}^\sharp$). $\hfill \square$\\

\section{The semi-simple filtration}

In this section every (convergent)  (a,b)-module is assumed to be regular.

\subsection{Semi-simple regular (a,b)-modules}

It is easy to see that a regular (a,b)-module $\mathcal{E}$ is {\bf simple}  (so by definition  has no non trivial {\em normal} sub-module) is either $\mathcal{E} = \{0\}$ or has  rank $1$.

\begin{defn}\label{ss 0}
Let \ $\mathcal{E}$ \ be a regular (a,b)-module. We say that \ $\mathcal{E}$ \ is {\bf semi-simple} if it is a sub-module of a finite direct sum of rank 1 regular (a,b)-modules.
\end{defn}

It is clear from this definition that a sub-module of a semi-simple (a,b)-module is semi-simple and that a (finite) direct sum of semi-simple (a,b)-modules is again semi-simple.

\parag{Remarks} \begin{enumerate}
\item We have already seen that a rank $1$ regular (a,b)-module is isomorphic to $ \mathcal{E}_\lambda := B[a]/B[a](a-\lambda b)$, and then a convergent (a,b)-module is simple if and only if its formal completion is simple. Then a semi-simple convergent (a,b)-module has a formal completion which is also semi-simple, that is to say which  is embeddable in a  finite direct sum of formal simple (a,b)-modules.
\item  A rather easy consequence of the classification of formal  rank \ $2$ \ (a,b)-modules given in \cite{[B.93]}  is that the rank \ $2$ \  (a,b)-modules defined in the \ $ \widehat{B} = \mathbb{C}[[b]]-$basis \ $x,y$ \  by the relations :
$$(a - (\alpha+p-1)b)x = y +  b^{p}y \quad {\rm and} \quad (a - \alpha b)y = 0$$
for any \ $\alpha\in \mathbb{C}$ \ and \ any \ $p \in \mathbb{N}^*$ \ are not semi-simple. We leave the verification of this point to the reader.\\
So the analog rank $2$ $B$-modules with $a$ defined by the same relations are also not semi-simple.
\end{enumerate}

\begin{lemma}\label{quot. ss 1}
Let \ $\mathcal{E}$ \ be an (a,b)-module which is direct sum of regular rank \ $1$ \ (a,b)-modules, and let \ $\mathcal{F} \subset \mathcal{E}$ \ be a rank  $1$ \ normal sub-module. Then \ $\mathcal{F}$ \ is a direct factor of \ $\mathcal{E}$ and we have $\mathcal{E} = \mathcal{F} \oplus \mathcal{H}$ where $\mathcal{H}$ is again a finite direct sum of regular rank \ $1$ \ (a,b)-modules.
\end{lemma}

\begin{cor}\label{quot. ss 2}
If \ $\mathcal{E}$ \ is a semi-simple regular (a,b)-module and \ $\mathcal{F}$ \ a normal sub-module of \ $\mathcal{E}$, the quotient \ $\mathcal{E}\big/\mathcal{F}$ \ is a (regular) semi-simple (a,b)-module.
\end{cor}

\parag{Proof of the lemma} Let  \ $\mathcal{E} = \oplus_{j=1}^k \ \mathcal{E}_{\alpha_j}$ \ and assume that \ $\mathcal{F} \simeq \mathcal{E}_{\beta}$. Let \ $e_j$ \ be a standard generator of \ $\mathcal{E}_{\alpha_j}$ \ and \ $e$ \ be  a standard generator of \ $\mathcal{E}_{\beta}$. Write $ e = \sum_{j=1}^k \ S_j(b)e_j $
and compute \ $(a - \beta b)e = 0$ \ using the fact that \ $(e_j, j \in [1,k])$, is a \ $B$-basis of \ $\mathcal{E}$ \ with  the relations \ $(a - \alpha_jb)e_j = 0$ \ for each \ $j$. We obtain, for each \ $j \in [1,k]$, the relation
$$ bS'_j - (\beta-\alpha_j)S_j = 0 .$$
So, if \ $\beta - \alpha_j$ \ is not in \ $\mathbb{N}$, we have \ $S_j = 0$. When \ $\beta = \alpha_j + p_j$ \ with \ $p_j \in \mathbb{N}$ \ we obtain \ $S_j(b) = \rho_jb^{p_j}$ \ for some \ $\rho_j \in \mathbb{C}$. As we assume that \ $e$ \ is not in \ $b\mathcal{E}$, there exists at least one \ $j_0 \in [1,k]$ \ such that \ $p_{j_0} = 0$ \ and \ $\rho_{j_0} \not= 0 $. In the case where there exists only one $j_0$ with $p_{j_0} = 0$ and $\rho_{j_0} \not= 0$  it is clear that we have
$$ \mathcal{E} = \mathcal{F} \oplus \big(\oplus_{j\not= j_0} \mathcal{E}_{\alpha_j}\big).$$
If there are many such $j_0$ then we are reduced to the case where $\mathcal{E}$ is the direct sum of several copies of  $ \mathcal{E}_\beta$ and where $e$ is a  (complex) linear combination of the standard generators. This case is also obvious, concluding the proof. $\hfill \blacksquare$

\parag{Proof of the corollary} We argue by induction on the rank of \ $\mathcal{F}$. In the rank \ $1$ \ case, we have \ $\mathcal{F} \subset \mathcal{E} \subset \mathcal{H} :=  \oplus_{j=1}^k \ \mathcal{E}_{\alpha_j}$. Let \ $\tilde{\mathcal{F}}$ \ the normalization of \ $\mathcal{F}$ \ in \ $\oplus_{j=1}^k \ \mathcal{E}_{\alpha_j}$. Then the lemma shows that there exists a \ $j_0 \in [1,k]$ \ such that 
 $$\mathcal{H} = \tilde{\mathcal{F}} \oplus\big( \oplus_{j\not= j_0} \mathcal{E}_{\alpha_j}\big).$$
 Then, as \ $\tilde{\mathcal{F}} \cap \mathcal{E} = \mathcal{F}$, the quotient map \ $\mathcal{E} \to \mathcal{H}\big/\tilde{\mathcal{F}}$  induces an injection of \ $\mathcal{E}\big/\mathcal{F}$ \ in a direct sum of regular rank \ $1$ \ (a,b)-modules. So \ $\mathcal{E}\big/\mathcal{F}$ \ is semi-simple.\\
 Assume now that the result is proved for \ $\mathcal{F}$ \ with rank \ $\leq d-1$ \ and assume that \ $\mathcal{F}$ \ has rank \ $d$. Then using a rank \ $1$ \ normal sub-module \ $\mathcal{G}$ \ in \ $\tilde{\mathcal{F}}$\footnote{Such a $\mathcal{G}$  exists as a consequence of  Corollary \ref{sol. Jordan}.}, the normalisation of $\mathcal{F}$ in $\mathcal{E}$, we obtain that \ $\mathcal{F}\big/(\mathcal{F}\cap \mathcal{G})$ \ is a  rank \ $d-1$ \ sub-module of \ $\mathcal{E}\big/\mathcal{G}$. Using the rank \ $1$ \ case we know that \ $\mathcal{E}\big/\mathcal{G}$ \ is semi-simple, and the induction hypothesis gives that 
  $$\mathcal{E}\big/\mathcal{F} = (\mathcal{E}\big/\mathcal{G})\Big/(\mathcal{F}\big/(\mathcal{F}\cap \mathcal{G}))$$
   is semi-simple. $\hfill \blacksquare$\\
   
   \begin{cor}\label{ss. 3} Let $\mathcal{E}$ be a simple pole semi-simple (a,b)-module. Then $\mathcal{E}$ is a direct sum of rank $1$ simple pole (a,b)-modules.
   \end{cor}
   
   \parag{Proof} This is obtained from Lemma \ref{quot. ss 1} by an easy induction on the rank of $\mathcal{E}$ which is left to the reader.$\hfill \blacksquare$

   \parag{Remark} If a regular (a,b)-module is semi-simple, its Bernstein polynomial has only simple roots.

 \subsection{The semi-simple filtration}

\begin{defn}\label{element ss 1}
Let \ $\mathcal{E}$ \ be a regular (a,b)-module and \ $x$ \ an element in \ $\mathcal{E}$. We shall say that \ $x$ \ is {\bf semi-simple} if \ $B[a]x$ \ is a semi-simple (a,b)-module.
\end{defn}

It is clear that any element in a semi-simple  (a,b)-module is semi-simple. The next lemma shows that the converse is true.

\begin{lemma}\label{element ss 2}
Let \ $\mathcal{E}$ \ be a regular (a,b)-module such that any \ $x \in \mathcal{E}$ \ is semi-simple. Then \ $\mathcal{E}$ \ is semi-simple. 
\end{lemma}

\parag{proof} Let \ $e_1, \dots, e_k$ \ be a \ $B$-basis of \ $\mathcal{E}$. Then each \ $B[a]e_j$ \ is semi-simple, and the map \ $\oplus_{j=1}^k B[a]e_j \to \mathcal{E} $ \ is surjective. So \ $\mathcal{E}$ \ is semi-simple thanks to  Corollary \ref{quot. ss 2} and the comment following Definition \ref{ss 0}. $\hfill \blacksquare$

\begin{lemma}\label{ss part 0}
Let \ $\mathcal{E}$ \ be a regular (a,b)-module. The subset \ $S_1(\mathcal{E})$ \ of semi-simple elements in \ $\mathcal{E}$ \ is a normal sub-module in \ $\mathcal{E}$.
\end{lemma}

\parag{proof} As a finite direct sum of semi-simple (a,b)-modules and also a quotient of a semi-simple (a,b)-module by a normal sub-module is semi-simple, it is clear that for  \ $x$ \ and \ $y$ \ semi-simple the sum \ $B[a]x + B[a]y$ \ is semi-simple. So \ $x + y$ \ is semi-simple. This implies that \ $S_1(\mathcal{E})$ \ is a sub-module of \ $\mathcal{E}$. \\
If \ $bx$ \ is in \ $S_1(\mathcal{E})$, then \ $B[a]bx$ \ is semi-simple. Then since $B[a]bx = bB[a]x $, $B[a]x $
  is also semi-simple, and then \ $S_1(\mathcal{E})$ \ is normal in \ $\mathcal{E}$. $\hfill \blacksquare$

\begin{defn}\label{ss part 1}
Let \ $\mathcal{E}$ \ be a regular (a,b)-module. Then the sub-module \ $S_1(\mathcal{E})$ \ of semi-simple elements in \ $\mathcal{E}$ \ will be called the {\bf\em semi-simple part} of \ $\mathcal{E}$.\\
Defining inductively  \ $S_j(\mathcal{E})$ \ as the pull-back in \ $\mathcal{E}$ \ of the semi-simple part of \ $\mathcal{E}\big/ S_{j-1}(\mathcal{E})$ \ for \ $j \geq 1$ \ with the initial condition \ $S_0(\mathcal{E}) = \{0\}$, we obtain an increasing  sequence of normal sub-modules in \ $\mathcal{E}$ \ such that \ $S_j(\mathcal{E})\big/S_{j-1}(\mathcal{E}) = S_1(\mathcal{E}\big/S_{j-1}(\mathcal{E}))$ is semi-simple. We shall call it the {\bf \em semi-simple filtration} of \ $\mathcal{E}$.\\
 The smallest integer \ $d$ \ such we have \ $S_d(\mathcal{E}) = \mathcal{E}$ \ will be called the {\bf\em nilpotent order} of \ $\mathcal{E}$ \ and we shall denote it \ $d(\mathcal{E})$.\\
\end{defn}

This filtration has the following properties
\begin{lemma} \label{basic filt.ss} 
Let $\mathcal{E}$ be a regular  (a,b)-module and let $\mathcal{F}$ be any sub-module. Then $S_j(\mathcal{F}) = S_j(\mathcal{E}) \cap \mathcal{F}$ for any $j \in \mathbb{N}^*$ with the convention that, for $j \geq d(\mathcal{E})+1$, we define $S_j(\mathcal{E}) = \mathcal{E}$.
In particular $d(\mathcal{F}) \leq d(\mathcal{E})$.\\
If $\mathcal{F}$ is normal  in $\mathcal{E}$, the quotient map $\mathcal{E} \to \mathcal{E}/\mathcal{F}$ sends $S_j(\mathcal{E})$ in $S_j(\mathcal{E}/\mathcal{F})$. This implies that $d(\mathcal{E}/\mathcal{F}) \leq d(\mathcal{E})$.\\
For any subset $\mathscr{A} \in \mathbb{C}/\mathbb{Z}$ and any $j$ we have $S_j(\mathcal{E})_{[\mathscr{A}]} = S_j(\mathcal{E}_{[\mathscr{A}]})$.
\end{lemma}

\parag{Proof} The inlusion $S_1(\mathcal{F}) \subset S_1(\mathcal{E}) $ is obvious, so $S_1(\mathcal{F}) = \mathcal{F}\cap S_1(\mathcal{E})$ holds true. \\
Assume that we have already proved that $S_j(\mathcal{F}) = S_j(\mathcal{E}) \cap \mathcal{F}$ for some $j \geq 1$. Then we have
\begin{align*}
&  S_{j+1}(\mathcal{F})/S_j(\mathcal{F}) = S_1(\mathcal{F}/S_j(\mathcal{F})) = S_1(\mathcal{F}/(S_j(\mathcal{E}) \cap \mathcal{F})) = \\
&  \big(\mathcal{F}/(S_j(\mathcal{E})\cap \mathcal{F})\big) \cap S_1(\mathcal{E}/S_j(\mathcal{E})) =\big( \mathcal{F}/(S_j(\mathcal{E})\cap \mathcal{F})\big) \cap (S_{j+1}(\mathcal{E})/S_j(\mathcal{E}))
\end{align*}
and this gives the equality $S_{j+1}(\mathcal{F}) = S_{j+1}(\mathcal{E}) \cap \mathcal{F}$ thanks to inductive assumption.\\
If $\mathcal{F}$ is normal the image of $S_1(\mathcal{E})$ by the quotient map $\mathcal{E} \to \mathcal{E}/\mathcal{F}$ is a  semi-simple sub-module in $\mathcal{E}/\mathcal{F}$ so is contained in $S_1(\mathcal{E}/\mathcal{F})$. \\
Assume that we have already proved that for each $h \in [1, j]$ the quotient map sends $S_h(\mathcal{E}) \to S_h(\mathcal{E}/\mathcal{F})$. Then it sends $S_j(\mathcal{E})/S_{j-1}(\mathcal{E})$ into $S_j(\mathcal{E}/\mathcal{F})/S_{j-1}(\mathcal{E}/\mathcal{F})$ and so
$S_1(\mathcal{E}/S_{j}(\mathcal{E})) = S_{j+1}(\mathcal{E})/S_j(\mathcal{E})) $ into $S_1(\mathcal{E}/\mathcal{F})/S_{j}(\mathcal{E}/\mathcal{F})) = S_{j+1}(\mathcal{E}/\mathcal{F})/S_j(\mathcal{E}/\mathcal{F})$. The conclusion follows.\\
Now if $d := d(\mathcal{E})$ we have $\mathcal{E} = S_d(\mathcal{E})$ and the quotient map sends $S_d(\mathcal{E}) = \mathcal{E}$ into $S_d(\mathcal{E}/\mathcal{F})$. Then $S_d(\mathcal{E}/\mathcal{F}) = \mathcal{E}/\mathcal{F}$.\\
To complete the proof, remark that the equalities, for any sub-module $\mathcal{F} $ of a regular (a,b)-module $\mathcal{E}$,  $S_j(\mathcal{E}) \cap \mathcal{F} = S_j(\mathcal{F})$ and $\mathcal{E}_{\mathscr{A}} \cap \mathcal{F} = \mathcal{F}_{\mathscr{A}}$ give
$$S_j(\mathcal{E}_{\mathscr{A}}) = S_j(\mathcal{E})\cap \mathcal{E}_{\mathscr{A}} = S_j(\mathcal{E})_{\mathscr{A}}. \qquad \qquad \qquad \qquad  \blacksquare$$

\parag{Remarks}
\begin{enumerate}
\item As \ $S_1(\mathcal{E})$ \ is the maximal semi-simple sub-module of \ $\mathcal{E}$ \ it contains any rank \ $1$ \ sub-module of \ $\mathcal{E}$. So \ $S_1(\mathcal{E}) = \{0\}$ \ happens if and only if \ $\mathcal{E} = \{0\}$.
\item Let $\mathcal{F}$ be a sub-module of  the regular (a,b)-module $\mathcal{E}$ such that $S_j(\mathcal{F}) = \mathcal{F}$. Then $\mathcal{F}$ is contained in $S_j(\mathcal{E})$ thanks to the previous lemma.
\item The semi-simple filtration of \ $\mathcal{E}$ \ is strictly increasing for \ $j \in [0, d(\mathcal{E})-1]$. 
\end{enumerate}

 The next lemma will help to compute the ranks of the various  $S_j(\mathcal{E})$ 
  
  \begin{lemma}\label{manque 1} Let $\mathcal{E}'$ be a sub-module of a regular (a,b)-module $\mathcal{E}$ such that $\mathcal{E}/\mathcal{E}'$ is a finite dimensional complex space. Then, for each $j$ the quotient $S_j(\mathcal{E})/S_j(\mathcal{E}')$ is also a finite dimensional complex vector space. So $S_j(\mathcal{E})$ and $S_j(\mathcal{E}')$ have the same rank on $B$ and $\mathcal{E}'$ and $\mathcal{E}$ have the same nilpotent order.
  \end{lemma}
  
  \parag{Proof} As we know that $S_j(\mathcal{E}') = S_j(\mathcal{E}) \cap \mathcal{E}'$ the quotient $S_j(\mathcal{E})/S_j(\mathcal{E}')$ is a sub-vector space of $\mathcal{E}/\mathcal{E}'$ and then has finite dimension. So  $S_j(\mathcal{E}')$ has finite co-dimension in $S_j(\mathcal{E})$ for each $j$ and then  has the same rank (as a $B$-module). \\
  Conversely, if $S_\delta(\mathcal{E}') = \mathcal{E}'$, the equality $S_\delta(\mathcal{E}') = S_\delta(\mathcal{E})\cap  \mathcal{E}'$ shows that $S_\delta(\mathcal{E})$ has  the same rank than $\mathcal{E}'$  and $\mathcal{E}$  and so must be equal to $\mathcal{E}$.$\hfill \blacksquare$\\

 \begin{cor}\label{primitive versus ss}
Let $\mathcal{E}$ be a regular (a,b)-module. If the  (a,b)-module $\mathcal{E}$ is semi-simple, then  $\mathcal{E}_{\mathscr{A}}$ is semi-simple for any subset $\mathscr{A} \subset \mathbb{C}/\mathbb{Z}$.\\
Conversely, if for any $\alpha \in \mathbb{C}/\mathbb{Z}$, $\mathcal{E}_{[\alpha]}$ is semi-simple, then $\mathcal{E}$ is semi-simple.
\end{cor}

\parag{Proof} The direct part is clear. \\
Since $\mathcal{E}$ is semi-simple if and only if $\mathcal{E}^\sharp$ is semi-simple, the converse follows from the finite co-dimension in $\mathcal{E}^\sharp$ of  the sub-module 
$\sum_{\alpha \in \mathscr{A}} \mathcal{E}_{[\alpha]}$ and the previous lemma.$\hfill \blacksquare$\\

  The next result shows that the ranks of the successive  quotients of the semi-simple filtration is non increasing.

  \begin{prop}\label{ss 2}
  Let $\mathcal{E}$ be a regular (a,b)-module and note $d := d(\mathcal{E})$ its nilpotent order. Then for each $j \in [1, d]$ we have
   $$rk(S_{j}(\mathcal{E})/S_{j-1}(\mathcal{E})) \geq rk(S_{j+1}(\mathcal{E})/S_j(\mathcal{E})).$$
  \end{prop}
  
  \parag{Proof} Remark that it is enough to prove the result when $\mathcal{E}$   has a simple pole, because 
  $S_j(\mathcal{E}^\sharp) = N_{\mathcal{E}^\sharp}(S_j(\mathcal{E})^\sharp) =  N_{\mathcal{E}^\sharp}(S_j(\mathcal{E})) $ for each $j$ so $S_j(\mathcal{E}^\sharp), S_j(\mathcal{E})^\sharp$ and $S_j(\mathcal{E})$ have the same rank.\\
  Note also that for $d(\mathcal{E}) = 1$ there is nothing to prove.\\
  We shall begin  by the case where $d(\mathcal{E}) = 2$. Then consider the exact sequence of (a,b)-modules  $0 \to S_1(\mathcal{E}) \to \mathcal{E} \to \mathcal{E}/S_1(\mathcal{E}) \to 0$ which give the exact sequence of finite dimensional vector spaces
  $$ 0 \to S_1(\mathcal{E})/bS_1(\mathcal{E}) \to \mathcal{E}/b\mathcal{E} \to \mathcal{E}/(S_1(\mathcal{E})+b\mathcal{E}) \to 0 $$
  which is compatible with the respective actions of $b^{-1}a$ on these quotients. But as $S_1(\mathcal{E})$ and $\mathcal{E}/S_1(\mathcal{E})$ are semi-simple, the action of $b^{-1}a$ on $S_1(\mathcal{E})/bS_1(\mathcal{E})$ and on $\mathcal{E}/(S_1(\mathcal{E})+b\mathcal{E})$ are semi-simple and the action of $b^{-1}a$ on $\mathcal{E}/b\mathcal{E}$ has a nilpotent part $N$ which satisfies $N^2 = 0$. So we have $ Im N \subset  Ker N$ and $\dim \mathcal{E}/b\mathcal{E} \leq 2\dim Ker N$. Since $Ker N = S_1(\mathcal{E})/bS_1(\mathcal{E})$\, because the (a,b)-module generated by a diagonal basis of $Ker N$ for $b^{-1}a$ generates a semi-simple sub-module of $\mathcal{E}$ (see Lemma \ref{element ss 2}), so it is equal to $S_1(\mathcal{E})$,  the conclusion follows for $d(\mathcal{E}) = 2$.\\
  To prove the general case, consider a geometric (a,b)-module $\mathcal{E}$ with $d(\mathcal{E}) \geq 3$ and an integer $j \in [2,d(\mathcal{E})-1]$. Then let $\mathcal{F} := S_{j+1}(\mathcal{E})/S_{j-1}(\mathcal{E})$. We have $d(\mathcal{F}) = 2$ and so
  $$ rk(\mathcal{F}) \leq 2 rk(S_1(\mathcal{F})) .$$
  Note  $s_h := rk(S_h(\mathcal{E})) $ for each $h $. The inequality proved  above  gives $s_j + s_{j+1} \leq 2 s_j$ because $S_1(\mathcal{F}) = S_j(\mathcal{E})/S_{j-1}(\mathcal{E})$ as $\mathcal{F} $ is a sub-module of $\mathcal{E}/S_{j-1}(\mathcal{E})$ which implies  $S_1(\mathcal{F}) = \mathcal{F} \cap S_1(\mathcal{E}/S_{j-1}(\mathcal{E})) = S_{j}(\mathcal{E})/S_{j-1}(\mathcal{E})$. This concludes the proof of the inequality $s_{j+1} \leq s_j$ for each $j$.$\hfill \blacksquare$\\
  
  The following corollary will be useful later on
  
  \begin{cor}\label{utile}
  Let $\mathcal{E}$ be a regular (a,b)-module which has a unique rank $1$ normal sub-module. Then, for any $j \geq 2$, $\mathcal{E}$ cannot have two distinct  normal rank $j$ sub-modules.
  \end{cor}
  
  \parag{Proof} 
   We begin by the case $j = 2$. We shall argue by contradiction. So let $G_1$ and $G_2$ be  two distinct normal sub-modules with rank $2$. By uniqueness of the normal rank $1$ sub-module $H$ of $\mathcal{E}$ we have $H \subset G_1 \cap G_2 $ since each $G_i$ must contain $H := S_1(\mathcal{E})$, and  moreover the intersection cannot be of rank $2$ since it is also  a normal sub-module. Then $G_1/H$ and $G_2/H$ are two distinct normal rank $1$ sub-modules of $S_2(\mathcal{E})/H$. Then $S_2(\mathcal{E}/S_1(\mathcal{E}))$ and then also $S_2(\mathcal{E})\big/S_1(\mathcal{E})$  have  rank at least two  while $S_1(\mathcal{E}) = H$ has rank $1$. This contradicts Proposition \ref{ss 2}.\\
Now we shall argue by induction on $j \geq 3$. So assume that $j \geq 2$ and the result proved for the rank  at most equal $j-1$. We shall argue again by contradiction. Then the quotient $\mathcal{E}/F_{j-2}$, where $F_{j-2}$ is the unique normal sub-module of rank $j-2$ in $\mathcal{E}$ thanks to our inductive assumption, has a unique rank $1$ normal sub-module and then has a unique normal rank $2$ normal submodule. If $G_1$ and $G_2$ are distinct normal sub-modules of rank $j$ in $\mathcal{E}$  then $G_1\cap G_2$ is equal to $F_{j-1}$, the unique normal rank $j-1$ sub-module of $\mathcal{E}$\footnote{Corollary \ref{sol. Jordan} implies that  any rank $j \geq 2$ normal sub-module of a regular (a,b)-module,  contains a rank $j-1$ normal sub-module.} and $F_{j-1}$ contains $F_{j-2}$. \\
  Then $G_i/F_{j-1}, i = 1,2$ are two distinct rank $1$ normal sub-modules of $\mathcal{E}/F_{j-1}$. Contradiction.  $\hfill \blacksquare$\\
  
  \parag{Remark} Note that under the hypothesis of the previous corollary the rank of $S_1(\mathcal{E})$ is $1$ (if $\mathcal{E} \not= \{0\}$). Then, the rank of $S_2(\mathcal{E})$ is $1$ or $2$. In the first case $S_2(\mathcal{E}) = S_1(\mathcal{E})$ and $\mathcal{E} = S_1(\mathcal{E})$ has rank $1$. If the rank of $S_2(\mathcal{E})$ is $2$, then the rank of $S_3(\mathcal{E})$ is $2$ or $3$. If the rank of $S_3(\mathcal{E})$ is $2$, then $S_2(\mathcal{E}) = \mathcal{E}$ and the rank of $\mathcal{E}$ is $2$. And so on $\dots$.
  So if $k$ is the rank of $\mathcal{E}$ we have $S_k(\mathcal{E}) = \mathcal{E}$ and $d(\mathcal{E}) = k$.\\
  Conversely, if $\mathcal{E}$ is regular with rank $k$ and such that $d(\mathcal{E}) = k$, the inclusions $S_j(\mathcal{E}) \subset S_{j+1}(\mathcal{E})$ are strict for $j \in [0,k-1]$ and each quotient $S_{j+1}(\mathcal{E})/S_j(\mathcal{E})$ has rank $1$. In particular $S_j(\mathcal{E})$ has rank $j$ and so $S_j(\mathcal{E})$  is the unique normal rank $j$ sub-module in $\mathcal{E}$. \\
  
 We leave the proof of the following corollary as an exercise for the reader :
  
  \begin{cor}\label{theme -2} 
  For any non zero regular (a,b)-module $\mathcal{E}$ we have the equivalence between the conditions:
  \begin{itemize}
  \item   $rank(S_1(E)) = 1$ and 
  \item $d(\mathcal{E}) = rank(\mathcal{E})$ 
  \item $\mathcal{E}$ is a $[\alpha]$-primitive theme.
  \item For each $j \in [0, rk(\mathcal{E})]$, $\mathcal{E}$  has a unique rank $j$ normal sub-module.$\hfill \blacksquare$\\
  \end{itemize}
  \end{cor}

  \section{Asymptotic expansions and geometric (a,b)-modules}
  
  \subsection{The Embedding Theorem}

Fix $\alpha \in ]0,1] \cap \mathbb{Q}$ and $N \in \mathbb{N}$. We shall consider convergent  asymptotic expansions of the type
\begin{equation}
\Xi_\alpha^{(N)} := \{ \sum_{j=0}^N\sum_{m = 0}^\infty  c_m^j s^{\alpha +m -1} \frac{(Log \,s)^j}{j!} \} 
\end{equation}
where we ask that exists $R > 0$ and a positive constant $C_R$ such that the following estimates hold 
\begin{equation} 
\vert c_m^j\vert \leq C_R R^m \quad \forall j \in [0,N]\quad \forall m \in \mathbb{N}.
\end{equation}
Of course, $\Xi_\alpha^{(N)}$ is a free   $\mathbb{C}\{s\}$-module with rank $(N+1)$ and  basis 
\begin{equation}
e_j := s^{\alpha-1}\frac{(Log \,s)^j}{j!}, \  j \in [0,N].
\end{equation}
We shall consider the dual Fr\'echet topology on $\Xi_\alpha^{(N)}$ deduced from the natural topology of $\mathbb{C}\{s\}$ (which is independent of the choice of the basis) and we define on $\Xi_\alpha^{(N)}$
the $\mathbb{C}$-linear continuous endomorphism\footnote{The continuity is an easy consequence of the formula 
$$ b(s^me_j) = \sum_{h=0}^j \frac{(-1)^{j-h}}{(\alpha+m)^{j-h+1}} s^{m+1}e_h$$ since $j \in [0,N], m \in \mathbb{N}$ and $\alpha > 0$.} defined inductively on $j \in [0,N]$ by the following rules:
\begin{enumerate}
\item For $j = 0$ and $m \in \mathbb{N} $ we put  $b(s^{\alpha + m-1}) := s^{\alpha + m}/(\alpha + m) $.
\item For $j \geq 1$ and $m \in \mathbb{N} $ we put
 \begin{equation*}
 b(s^{m}e_j) = \frac{s^{m+1}}{\alpha + m} e_j - \frac{b(s^m e_{j-1})}{\alpha + m}. \tag{E}
 \end{equation*}
\end{enumerate}
Of course this endomorphism corresponds to the  term-wise  integration (without constant) of these series.

Remind that we denote  $A := \mathbb{C}\{a\}$ and $B := \mathbb{C}\{\{b\}\}$ the closed sub-algebras of the continuous $\mathbb{C}$-linear endomorphism of 
$\mathbb{C}\{s\}$ generated by $a = \times s$ and $b := \int_0^s$.

\begin{prop}\label{expansion 1}
The $\mathbb{C}[b]$-action on $\Xi_\alpha^{(N)}$ defined above extends to a $B$-action and $\Xi_\alpha^{(N)}$ is a free $B$-module of rank $N+1$ with basis $e_j, j \in [0,N]$. Moreover, we have
$a\Xi_\alpha^{(N)} = b\Xi_\alpha^{(N)}$, and the action of $b^{-1}a$ on $\Xi_\alpha^{(N)}$ may be written as
$$ b^{-1}a =  \Delta + \mathcal{N} \quad  {\rm with} \quad  \mathcal{N}(e_0) = 0  \quad  {\rm and}   \quad \mathcal{N}(e_j) = e_{j-1} \quad  \forall j \in [1, N]$$
where $\mathcal{N}$ is $A$-linear and $B$-linear, where $\Delta(s^m e_j) = (\alpha + m)s^m e_j, \  \forall m \  \forall j$  and where $b^{-1}a,\, \Delta$ and $\mathcal{N}$  commute.
\end{prop}

\parag{Proof} The formula $(E)$ gives for each integers $m$ and $j$  
$$ as^me_j = (\alpha + m)bs^me_j + bs^me_{j-1} $$
so $\Delta(s^me_j) = (\alpha + m)s^me_j + \mathcal{N}(s^me_j)$ if we define the $A$-linear map $\mathcal{N}$ by putting $\mathcal{N}(e_j) = e_{j-1}$ with the convention $\mathcal{N}(e_0) = 0$.\\
Then it is clear that $\mathcal{N}$ satisfies $\mathcal{N}^{N+1} = 0$. The fact that $\Delta$ is $\mathbb{C}$-linear continuous and bijective gives $a\Xi_\alpha^{(N)} = b\Xi_\alpha^{(N)}$.\\
We want to check now that $\mathcal{N}$ is $b$-linear.\\
 Then  first remark that we have  $b^{-1}ab - bb^{-1}a = b^{-1}b(a+b) -a = b$. \\
 We shall prove by induction on $j \geq 0$ that $\Delta b - b\Delta = b$.\\
For $j = 0$ we have $$(\alpha + m)\Delta b s^me_0 = \Delta(s^{m+1}e_0) = (\alpha + m +1)s^{m+1}e_0 $$
 and 
 $$(\alpha + m)b\Delta s^me_0 = (\alpha + m)^2bs^me_0 = (\alpha + m)s^{m+1}e_0$$
   and so 
$$(\alpha + m)(\Delta b - b\Delta)(s^me_0) = s^{m+1}e_0 = (\alpha+m)b(s^me_0).$$
Assume now $j \geq 1$ and the relation $\Delta b - b\Delta = b$ proved for $s^me_{j-1}$.\\
 Then we have, using the induction hypothesis:
\begin{align*}
& \Delta(bs^me_j) = \frac{1}{\alpha+m}\Delta\big( s^{m+1}e_j - b(s^me_{j-1}\big) = \frac{\alpha+m+1}{\alpha+m}s^{m+1}e_j - \frac{1}{\alpha+m}\Delta(bs^me_{j-1}) \\
&  \Delta(bs^me_j) = \frac{\alpha+m+1}{\alpha+m}s^{m+1}e_j  - bs^me_{j-1} - \frac{1}{\alpha+m}bs^me_{j-1} \quad {\rm and \ also} \\
& b\Delta(s^me_j) = b (\alpha+m)s^me_j = s^{m+1}e_j - bs^me_{j-1} \quad {\rm so \ we \ obtain} \\
&( \Delta b - b\Delta)(s^me_j) = \frac{1}{\alpha+m}\big( s^{m+1}e_j - bs^me_{j-1}\big) = bs^me_j.
\end{align*}
This gives the $b$-linearity of $\mathcal{N} = b^{-1}a - \Delta$, and then the $B$-linearity by continuity.\\
Note that $[b^{-1}a, b] = a$ and the $A$-linearity of $\mathcal{N}$ gives $[\Delta, a] = a$ which is easy to check directly.\\
Using the fact that $B$ acts on $\Xi_\alpha^{(N)}$\footnote{The easy estimates corresponding to this assertion is left to the reader.} and that  $b$ is injective on $\Xi_\alpha^{(N)}$, we conclude, as $(e_j, j \in [0,N])$ induces a basis of the vector space $\Xi_\alpha^{(N)} /b\Xi_\alpha^{(N)}  = \Xi_\alpha^{(N)} /a\Xi_\alpha^{(N)} $, and then  $\Xi_\alpha^{(N)} $ is a free $B$-module with basis $(e_j, j \in [0,N])$.$\hfill \blacksquare$\\

Let $\mathscr{A}$ be a finite subset of $]0,1] \cap \mathbb{Q}$ and let $V$ be a finite dimensional complex space, we define the free finite type $B$-module
\begin{equation}
 \Xi_{\mathscr{A}}^{(N)} := \oplus_{\alpha \in \mathscr{A}} \Xi_\alpha^{(N)} 
 \end{equation}
 which is also a free finite type $A$-module. Then, defining the action of $B$ and of $A$ as the identity on $V$, the tensor product $ \Xi_{\mathscr{A}}^{(N)} \otimes_{\mathbb{C}} V$ is again a free finite type $B$-module which is also a free finite type $A$-module. Of course this implies that it is a (simple pole) convergent (a,b)-module.
 
 \begin{defn}\label{9/2}
 For $\alpha \in ]0, 1[ \cap \mathbb{Q}$ we define $S\Xi_{\alpha}^{(N)} = \Xi_{\alpha}^{(N)}$ and for $\alpha = 1$ we put  $S\Xi_1^{(N)} := \Xi_1^{(N+1)} \big/\Xi_1^{(0)} \simeq  \Xi_1^{(N)}$.\\
 Then for a finite subset $\mathcal{A}$ in $]0, 1] \cap \mathbb{Q}$ and $V$ a finite dimensional complex space we define
 $$ S\Xi_{\mathcal{A}}^{(N)}\otimes V := \oplus_{\alpha \in \mathcal{A}} \  S\Xi_{\alpha}^{(N)}\otimes V .$$
 \end{defn}
 
  Note that, in fact $S\Xi_{\alpha}^{(N)}$ is always isomorphic to $\Xi_{\alpha}^{(N)}$, so the previous considerations remains true with this small change which corresponds to the fact that we are only interested by the "singular part" in the asymptotic expansions. 
 
 \begin{defn}\label{geom. 1}
 We say that a left $B[a]$ sub-module $\mathcal{E}$  of  \ $S\Xi_{\mathscr{A}}^{(N)} \otimes_{\mathbb{C}} V$  for some $N \in \mathbb{N}$ and some finite dimensional complex vector space $V$ is a {\bf geometric} (a,b)-module.\\
  A geometric (a,b)-module $\mathcal{E}$  of the form $B[a]\varphi$ for some $\varphi$ in some $S\Xi^{(N)}_{\mathscr{A}}\otimes_\mathbb{C} V $ is  called a {\bf fresco}.\\
A fresco $\mathcal{E}$  will be called a {\bf theme} when it may be written $\mathcal{E} = B[a]\varphi$ with $\varphi \in S\Xi^{(N)}_{\mathscr{A}}$ (so with $V = \mathbb{C}$).
 \end{defn}
 
 \parag{Remarks}
  \begin{enumerate}
  \item Of course a geometric (a,b)-module is a regular (a,b)-module because it is a free finite rank $B$-module, stable by the action of $a$ which is continuous, and which is regular since any  $S\Xi_{\mathscr{A}}^{(N)} \otimes_{\mathbb{C}} V$ is a simple pole (a,b)-module.
  \item The saturation of a geometric (a,b)-module $\mathcal{E}$ is again a geometric (a,b)-module since the stability of  any  \ $S\Xi_{\mathscr{A}}^{(N)} \otimes_{\mathbb{C}} V$  by $b^{-1}a$ implies  the inclusion of  $\mathcal{E}^\sharp$  in $S\Xi_{\mathscr{A}}^{(N)} \otimes_{\mathbb{C}} V$  when  $\mathcal{E}$ is in $S\Xi_{\mathscr{A}}^{(N)} \otimes_{\mathbb{C}} V$.
  \item It is easy to see that the Bernstein polynomial of $S\Xi_\alpha^{(N)}$ is equal to $(x +\alpha)^{N+1}$ and then the roots of the  Bernstein polynomial of a geometric (a,b)-module are negative rational numbers. 
    \item For $\alpha = -1$ (which is exclude in our definitions above) the free  rank $1$ $A$-module $(1/s)\mathbb{C}\{s\}$ is not stable by $b$, but the free rank $2$ $A$-module 
    $$(1/s)\mathbb{C}\{s\} \oplus (Log\, s)\mathbb{C}\{s\}$$
     is stable by $b$ if we defined
    $b(1/s) := Log\, s$ and $b(Log\, s) := sLog\, s - s$ which gives $$ba(1/s) = ab(1/s) -b^2(1/s) = sLog\, s -(sLog\, s - s) = s = b(1).$$
    So it coincides with the free rank $2$ $B$-module $\mathcal{E}$ with basis $1/s$ and $1$. \\
    But  its saturation by $b^{-1}a$ inside $\mathcal{E} \otimes_B B[b^{-1}]$ contains the non zero element $b^{-1}a(1/s) = b^{-1}[1]$ which satisfies $ab^{-1}[1]= 0$. So $a$ is not injective on $\mathcal{E}^\sharp$ which is not a free $A$-module although it is a free $B$-module.
  \end{enumerate}
 
\begin{thm}{\rm \bf [The Embedding Theorem]}\label{Emb. Th.}\label{ET}
 Let $\mathcal{E}$ be a regular  (a,b)-module such that its Bernstein polynomial has negative rational roots. Then there exists a finite subset $\mathscr{A} $  in $\mathbb{Q} \cap ]0,1]$,  a  finite dimensional complex vector space $V$, an integer $N \in \mathbb{N}$ and a $B[a]$-linear injective map 
 $f : \mathcal{E} \to S\Xi^{(N)}_\mathscr{A} \otimes V$.\\
 \end{thm}
 
 \parag{Important remark} An obvious consequence of the previous theorem is  the  following equivalent definition of a geometric (a,b)-module:
 \begin{itemize}
 \item A regular (a,b)-module is {\bf geometric} if and only if the roots of its Bernstein polynomial are negative rational numbers (compare with \cite{[M.74]} and  \cite{[K.76]}).
 \end{itemize}
 
 In the section 5.2  we describe, for a given geometric (a,b)-module $\mathcal{E}$, what are the smallest $N, \dim V$ and subset  $\mathscr{A}$ for which there exists an embedding of $\mathcal{E}$ in $S\Xi_\mathscr{A}^{(N)} \otimes_\mathbb{C} V$.\\
 
 The proof of the Embedding  Theorem will use the following lemmas\footnote{This easy  lemma is missing in the proof of Theorem 4.2.12 in \cite{[B.09]}.}.
 
  \begin{lemma}\label{manquant 1}
 Let $\gamma$ be a positive rational number. Then for any element $y$ in $ S\Xi^{(N)}_\mathscr{A}\otimes V$ there exists $x \in  S\Xi^{(N+1)}_\mathscr{A}\otimes V$ such that $ (a - \gamma b)x = by$.
 \end{lemma} 
 
 \parag{Proof} It is enough to consider the case $V = \mathbb{C}$ and $\mathscr{A} := \{\alpha\}$. \\ 
 Assume first that $\gamma \not\in \alpha + \mathbb{N}$.
 Since $\Delta - \gamma : S\Xi_\alpha^{(N)} \to S\Xi_\alpha^{(N)}$ is bijective and $\mathcal{N}^{N+1} = 0$ the map  $\Delta - \gamma - \mathcal{N}$ is also bijective and the formula
 $ a - \gamma b = b(\Delta - \gamma -\mathcal{N})$ implies the result, since $\Delta$ and $\mathcal{N}$ commute.\\
 If $\gamma = \alpha + m_0$ for some $m_0 \in \mathbb{N}$,  let $Z$ be the closed $\mathbb{C}$-linear span of the vectors $s^me_j$ for $ m \in \mathbb{N} \setminus \{m_0\}$ and $j \in [0,N]$. Then $\Delta -\gamma$ is continuous bijective on $Z$ and $\mathcal{N}$  satisfies $\mathcal{N}(Z) \subset Z$ and $\mathcal{N}^{N+1} = 0$. Now the formula  $ a - \gamma b = b(\Delta - \gamma -\mathcal{N})$ implies the equality $(a - (\alpha + m_0)b)Z = b(\Delta -\gamma - \mathcal{N})Z$. So, to complete the proof, it is enough to shows that $(a - (\alpha + m_0)b)(\Xi^{(N+1)}_\alpha)$ contains the vectors $s^{m_0+1}e_j, j \in [0,N]$.\\
  This is given by the formulas
   $$(a - (\alpha+m_0)b)(s^{m_0}e_j) = b(s^{m_0}e_{j-1}) $$
  and an induction on  $j$ in $[1,N+1]$.$\hfill \blacksquare$\\

  \parag{Proof of the Embedding Theorem} Remark first that it is enough to prove the existence of an embedding in the case where $\mathcal{E}$ has  simple pole.\\
 We  shall make an induction on the rank of $\mathcal{E}$ assumed to have a simple pole.\\
   In rank $1$ we have $\mathcal{E} \simeq \mathcal{E}_\alpha$ with $\alpha \in  \mathbb{Q}^{+*}$ and  since $\mathcal{E}_\alpha \simeq B[a] s^{\alpha-1}$ is a sub-module of $ S\Xi^{(0)}_{[\alpha]} $, this gives the desired embedding, where $[\alpha]$ is the class of $\alpha$ in  $\mathbb{Q}/\mathbb{Z} \simeq  \mathbb{Q} \cap ]0,1]$. \\
 Assume that the existence of such an embedding has been proved when $\mathcal{E}$ has rank $(k-1)$ with $k \geq 2$. Then consider a rank $k$ simple pole (a,b)-module $\mathcal{E}$ and let $\mathcal{F}$ be a rank $(k-1)$ normal sub-module in $\mathcal{E}$ (the existence of such an $\mathcal{F}$ is consequence of the existence of J-H sequence for $\mathcal{E}$;  see \cite{[B.93]}. The convergent case is analogous; see  Corollary \ref{sol. Jordan}). Then $\mathcal{E}/\mathcal{F}$ is isomorphic to $\mathcal{E}_\alpha$ for some $\alpha \in \mathbb{Q}^{+*}$ because it has rank $1$ and the only root $-\alpha$  of its  Bernstein polynomial is a root of the Bernstein polynomial of $\mathcal{E}$ (see Remark 3 following Definition \ref{Bernstein 0}). Let $e$ be an element in $\mathcal{E}$ such that its image in $\mathcal{E}/\mathcal{F}$ is the standard generator\footnote{ a non zero element $\varepsilon$  in $\mathcal{E}_\alpha$ such that $a\varepsilon = \alpha b\varepsilon$.} of $\mathcal{E}_\alpha$. Then, since $\mathcal{E}$ has a simple pole and $\mathcal{F}$ is normal, we have $z := (a - \alpha b)e$ is in $\mathcal{F}\cap b\mathcal{E} = b\mathcal{F}$. Then we may write $z = by$ with $y \in \mathcal{F}$.\\ 
 Now our induction hypothesis gives us an injective $B[a]$-linear embedding 
 $$g: \mathcal{F} \to S\Xi^{(N)}_\mathscr{A} \otimes V.$$
  Then $g(z)$ is in $bS\Xi^{(N)}_\mathscr{A} \otimes V$ and the lemma \ref{manquant 1} gives an $x \in S\Xi^{(N+1)}_\mathscr{A} \otimes V$ such that $(a - \alpha b)x = g(z)$. \\
   Define  a $B$-linear map  $f : \mathcal{E} \to S\Xi^{(N+1)}_\mathscr{A} \otimes (V\oplus \mathbb{C} \varepsilon)$ by
 \begin{align*} 
 &  f(t) = g(t)  \quad {\rm  when} \quad  t \in \mathcal{F} \quad  {\rm  and \ by} \\
 &  f(e) = x \oplus (s^{\alpha-1}\otimes \varepsilon), \quad  (resp. \quad  x \oplus (s^{\alpha-1}Log\, s)\otimes \varepsilon) \ {\rm for} \  \alpha \in \mathbb{N}^*),
 \end{align*}
 using the direct sum decomposition (as $B$-module) $\mathcal{E} = \mathcal{F} \oplus Be$.\\
 We shall verify the $a$-linearity of $f$ and also its injectivity.\\
 Let $\sigma$ be in $B$ and $t$ in $\mathcal{F}$. We have $a(t \oplus \sigma e) = (at + \sigma(a - \alpha b)e) \oplus \alpha\sigma be$. Then
 \begin{align*}
 & f(a(t\oplus \sigma e)) = g(at + \sigma(a - \alpha b)e) + \alpha\sigma b(x + s^{\alpha-1}\otimes \varepsilon) \quad (resp. \quad  x \oplus (s^{\alpha-1}Log\, s \otimes \varepsilon))  \\
 & \qquad = ag(t) + \sigma g(z) + \sigma s^\alpha \otimes \varepsilon - \sigma(a - \alpha b)x + \sigma ax \quad (resp. \quad  \sigma (s^{\alpha}Log\, s \otimes \varepsilon)) \\
 & \qquad = a f(t \oplus \sigma e) \quad {\rm since} \quad g(z) = (a -\alpha b)(x).
 \end{align*} 
 So the $a$-linearity is proved.\\
 As $g$ is injective, $f(t \oplus \sigma e) = 0$ implies $t +  \sigma x = 0$ and  $$ \sigma s^{\alpha-1}\otimes \varepsilon = 0 \quad (resp.\quad \sigma s^{\alpha-1}Log\, s \otimes  \varepsilon = 0 \quad {\rm when} \quad \alpha \in \mathbb{N}^* ).$$
 But since $W := V \oplus \mathbb{C}\varepsilon$, this implies $\sigma = 0$ and then $t = 0$, concluding the existence of an embedding for $\mathcal{E}$.$\hfill \blacksquare$\\
 
 \begin{cor}\label{theme -1}
 A $[\alpha]$-primitive theme of rank $\geq 2$ is not semi-simple. In fact its semi-simple part has rank $1$.
 \end{cor}
 
 \parag{Proof} Let $\mathcal{E} := B[a]\varphi \subset S\Xi_{[\alpha]}^{(N)}$ a rank $k$  $[\alpha]$-primitive theme which is semi-simple. Then $\mathcal{E}^\sharp$ is isomorphic to a direct sum of $\mathcal{E}_\beta$ for a finite set of $\beta$ in $\alpha + \mathbb{N}$ (may be with repetitions), thanks to Proposition \ref{ss. 3}. Since  $S\Xi_{[\alpha]}^{(N)}$ has a simple pole, the inclusion of $\mathcal{E}$ in $S\Xi_{[\alpha]}^{(N)}$ extends to a $B[a]$-linear map $j^\sharp : \mathcal{E}^\sharp \to S\Xi_{[\alpha]}^{(N)}$. But for any $\beta \in \alpha + \mathbb{N}$ a $B[a]$-linear map $\mathcal{E}_\beta \to S\Xi_{[\alpha]}^{(N)}$ has its image in $S\Xi_{[\alpha]}^{(0)}$ because $\mathbb{C}s^\beta$ (resp. $\mathbb{C} s^\beta Log\, s$ for $\alpha = 1$) is the vector space of solutions of the equation $(a - \beta b)x = 0$ for $x \in S\Xi_{[\alpha]}^{(N)}$. Then the image of $j^\sharp$ has rank at most $1$ and $k \leq 1$.$\hfill \blacksquare$\\
 
 \begin{cor}\label{theme -2}
 Let $\mathcal{E}$ be an $[\alpha]$-primitive theme of rank $k$. Then $d(\mathcal{E}) = k$.
 \end{cor}
 
 \parag{Proof} Remark first that we have an isomorphism $S\Xi_{[\alpha]}^{(N)}\big/S\Xi_{[\alpha]}^{(N-p)} \simeq S\Xi_{[\alpha]}^{(p-1)}$ which is given by $e_j \mapsto e_{j-p}$ for $j \in [N-p,N]$ (and $e_j \mapsto 0$ for $j \in [0,N-p]$), where $e_j := s^{\alpha-1}(Log\, s)^j/j!$ for $j \in [0,N]$ for $\alpha \not= 1$  and $e_j = Log\, s)^{j+1}\big/(j+1)!$ for $\alpha = 1$ and  for $j \in [0,N]$.  \\
Let $\mathcal{E}$ be an $[\alpha]$-primitive theme of rank $k\geq 2$. Then $\mathcal{E}/S_1(\mathcal{E})$ is an $[\alpha]$-primitive theme of rank $k-1$ since the proof of the previous corollary shows that an embedding of $\mathcal{E}$  in $S\Xi_{[\alpha]}^{(N)}$ gives an embedding of $\mathcal{E}/S_1(\mathcal{E})$ in $S\Xi_{[\alpha]}^{(N)}\big/S\Xi_{[\alpha]}^{(0)}\simeq S\Xi_{[\alpha]}^{(N-1)}$. Then 
 $$ \mathcal{E}/S_1(\mathcal{E})\big/S_1(\mathcal{E}/S_1(\mathcal{E})) \simeq \mathcal{E}/S_2(\mathcal{E}) $$
 is a rank $k-2$ $[\alpha]$-primitive theme. We obtain, continuing in this way, that for any $j \in [0, k-1]$ that $\mathcal{E}/S_j(\mathcal{E})$ is a rank $k-j$ $[\alpha]$-primitive theme and then that 
 $\mathcal{E} = S_k(\mathcal{E})$ and $d(\mathcal{E}) = k$.$\hfill \blacksquare$

 \subsection{Complements to the Embedding Theorem}
 
  We give now some complements to the Embedding Theorem \ref{Emb. Th.}.
 
 \begin{prop}\label{cas Asympt.} 
 Let $\mathscr{A}$ be a finite subset in $\mathbb{Q} \cap ]0,1]$, $V$ be a non zero finite dimensional vector space and $N$ be a non negative integer.\\
  Then the semi-simple filtration of $S\Xi^{(N)}_\mathscr{A} \otimes V$ is given by the equality
   $$S_j(S\Xi^{(N)}_\mathscr{A} \otimes V) = S\Xi^{(j-1)}_\mathscr{A} \otimes V$$
    for each $j \in [1, N+1]$, and the nilpotent order of this geometric (a,b)-module is $N+1$.
 \end{prop}
 
The proof of this proposition will use the following  three  lemmas.
 
  \begin{lemma}\label{theme 1}
  Fix  $\alpha \in \mathbb{Q} \cap ]0,1]$ and $ \beta := \alpha + m$ with  $m \in \mathbb{N}$. Let $P$ be a degree $M$ polynomial. Then there exists a unique degree $M+1$ polynomial $Q$
  without constant term such that $(a - \beta b)(s^{\beta-1}Q(Log\, s)) = s^\beta P(Log\, s)$.
  \end{lemma} 
  
  \parag{Proof} An elementary computation shows that $Q$ is the primitive vanishing at $0$ of the polynomial $\beta P + P'$ which has degree $M$ since $\beta > 0$. So $Q$ is a polynomial of degree $M+1$.$\hfill \blacksquare$ \\
  
  \begin{lemma}\label{theme 2}
We keep the notations of the previous lemma and assume $M \geq 1$.   Let $\varphi := s^{\beta-1}(Log\, s)^M + \psi$ where $\psi$ is in $\Xi^{(M-1)}_{[\alpha]}$. \\
Then the degree in $Log\, s$ of $(a - \gamma b)\varphi$ is equal to $M$ for $\gamma \not= \beta$ and $M-1$ for $\gamma = \beta$.
\end{lemma}

\parag{Proof} For $\gamma \not= \beta$ we have 
$$(a - \gamma b)s^{\beta-1}(Log\, s)^M = (1 - \gamma/\beta) s^\beta(Log\,s)^M \ {\rm modulo} \ \Xi^{(M-1)}_{[\alpha]} $$
 and the conclusion follows.\\
For $\gamma = \beta$ we have
 $$(a - \beta b) \varphi = Mb(s^{\beta-1}(Log\, s)^{M-1}) + (a-\beta b)\Xi^{(M-1)}_{[\alpha]}.$$
  But as $M \geq 1$ the Lemma \ref{theme 1} implies that $s^{\beta}(Log\, s)^{M-1}$ is not in
$(a-\beta b)\Xi^{(M-1)}_{[\alpha]}$ since the kernel of $a - \beta b$ is $\mathbb{C} s^{\beta-1}$. Now the equality
 $$b(s^{\beta-1}(Log\, s)^{M-1} = \frac{1}{\beta} s^\beta(Log\, s)^{M-1} - \frac{M-1}{\beta}b(s^{\beta-1}(Log\, s)^{M-2}) $$
 shows that the term $s^\beta(Log\, s)^{M-1}$ has still a non zero coefficient in $(a - \beta b)\varphi$ concluding the proof.$\hfill \blacksquare$\\

  \begin{lemma} \label{theme 3}
 Fix  $\alpha \in \mathbb{Q} \cap ]0,1]$ and $M \in \mathbb{N}$. Let $\varphi$ be in $S\Xi^{(M)}_{[\alpha]} $ of the form 
 $$ \varphi = s^{\beta-1}(Log\, s)^M + \psi \quad {\rm for} \ \alpha \not= 1 \quad  {\rm or} \quad s^{\beta-1}(Log\, s)^{M+1} + \psi  \quad {\rm for} \ \alpha = 1$$
 with $\beta $ in $\alpha + \mathbb{N}$ and $\psi \in S\Xi^{(M-1)}_{[\alpha]}$. Then $B[a]\varphi$ is a $[\alpha]$-primitive theme of rank $M+1$.
 \end{lemma}
 
 \parag{Proof} By definition $B[a]\varphi$ is a $[\alpha]$-primitive theme. We shall argue by induction on $M \geq 0$. For $M = 0$ since $S\Xi_{[\alpha]}^{(0)} \simeq \mathcal{E}_\alpha$ is rank $1$, the point is to prove that 
 $B[a]\varphi$ is not a finite dimensional complex space. As $a$ is injective this is clear.\\
 Assume that $M \geq 1$ and that the lemma is proved for $M-1$.  Thanks to  Lemma  \ref{theme 2} and the induction hypothesis we know that $B[a](a -\beta b)\varphi$ is a rank $M$ theme. Then the exact sequence
  $$0 \to  B[a](a -\beta b)\varphi \to B[a]\varphi \to Q \to 0$$
   where $Q$ is a quotient of $\mathcal{E}_\beta$, shows that the rank of $B[a]\varphi$ is either $M+1$ or $M$. If the rank is $M$, then, thanks to our induction hypothesis,  there exists $m \in \mathbb{N}$ such that $b^m\varphi$ is in $B[a](a-\beta b)\varphi$ which is contained in $S\Xi_{[\alpha]}^{(M-1)}$. This is clearly impossible. Then the rank of $B[a]\varphi $  is $M+1$.$\hfill \blacksquare$\\
   
   \parag{Remark} If, in the situation of the previous lemma,  $\varphi = S(b)s^{\beta-1}(Log\, s)^M + \psi $ for some $S$ invertible in $B$, then $\varphi$ also generates a $[\alpha]$-primitive theme.\\

 \parag{Proof of Proposition \ref{cas Asympt.}}  As $S\Xi^{(0)}_\mathscr{A}  \simeq \oplus_{\alpha \in \mathscr{A}} \mathcal{E}_\alpha$ it is clear that $S\Xi^{(0)}_\mathscr{A} \otimes V$ is semi-simple as  a finite direct sum of such $\mathcal{E}_\alpha$. So it is contained in $S_1(\mathcal{E})$.\\
  Conversely, if $\varphi \in S\Xi^{(N)}_\mathscr{A} \otimes V$ has degree at least equal to $1$ in $Log \, s$ (we deduce $1$ to the degree in $Log\, s$ when $\alpha = 1$), then Lemma \ref{theme 3} implies (by choosing a convenient linear form on $V$) that  some $[\alpha]$-primitive  sub-module $B[a]\varphi$ has a  $[\alpha]$-primitive quotient theme of rank $\geq 2$ for some $[\alpha] \in \mathscr{A}$, so is not semi-simple (see Corollary \ref{theme -1}). So $S_1(\mathcal{E}) = S\Xi^{(0)}_\mathscr{A} \otimes V$.\\
  Assume that $j \geq 1$ and that we have proved the equality  $$S_j(S\Xi^{(N)}_\mathscr{A} \otimes V) = S\Xi^{(j-1)}_\mathscr{A} \otimes V.$$
  Now, as $S\Xi^{(j)}_\mathscr{A} \otimes V\big/ S\Xi^{(j-1)}_\mathscr{A} \otimes V$ is semi-simple, since it is isomorphic to $S\Xi^{(0)}_\mathscr{A} \otimes V$, we obtain that $S\Xi^{(j)}_\mathscr{A} \otimes V \subset S_{j+1}(S\Xi^{(N)}_\mathscr{A} \otimes V)$ by the definition of $S_{j+1}(\mathcal{E})$. To complete our induction it is enough to prove that if $\varphi \in S\Xi^{(N)}_\mathscr{A} \otimes V$ has degree $j+1$ in $Log\, s$ (degree $j+2$ for $\alpha = 1$)  then $\varphi$ is not in $S_{j+1}(S\Xi^{(N)}_\mathscr{A} \otimes V)$. But under this hypothesis $B[a]\varphi$, thanks to Lemma  \ref{theme 3}, admits a quotient which is a rank $j+2$ theme.

  Then thanks to Corollary \ref{theme -2} $d(B[a]\varphi) =  j+2$ and  $\varphi$ is not in $S_{j+1}(S\Xi^{(N)}_\mathscr{A} \otimes V)$.$\hfill \blacksquare$\\

   As a consequence of the previous proposition, using Lemma \ref{basic filt.ss}, we obtain that  for any sub-module  $\mathcal{E} \subset (S\Xi^{(N)}_\alpha\otimes V)$ the equality  $S_j(\mathcal{E}) = \mathcal{E} \cap (S\Xi^{(j-1)}_{[\alpha]}\otimes V)$ for each $j \in [1, d(\mathcal{E})]$. \\

  The next proposition is also a complement to the Embedding Theorem \ref{Emb. Th.}.
  
  \begin{prop} \label{Manquant 2}
Let $\mathcal{E}$ be a geometric (a,b)-module and assume that $S_1(\mathcal{E})$ may be embedded in $S\Xi^{(0)}_{\mathscr{A}}\otimes V$. Then  we can extend  this embedding to an embedding of $\mathcal{E}$  in $S\Xi^{(N)}_{\mathscr{A}}\otimes V$ with $N := d(\mathcal{E}) - 1$.
\end{prop} 

\parag{Proof} To simplify the notation, we shall write $d := d(\mathcal{E})$ the nilpotent order of $\mathcal{E}$. Remark that it is enough to prove the result when $\mathcal{E}$ has a simple pole because any embedding of a geometric (a,b)-module in some $S\Xi^{(N)}_{\mathscr{A}}\otimes V$  extends to an embedding of  $\mathcal{E}^\sharp$, as any $S\Xi^{(N)}_{\mathscr{A}}\otimes V$ has a simple pole.\\
 Now as a geometric  simple pole (a,b)-module decomposes as a direct sum of its $[\alpha]$-primitive parts when $\alpha$ describes $\mathbb{Q} \cap ]0, 1]$, we may assume that $\mathcal{E}$ is $[\alpha]$-primitive.\\
   We shall prove the result by induction on $d$, assuming that $\mathcal{E}$ has a simple pole and  is $[\alpha]$-primitive (see Theorem \ref{DT}). The case $d= 1$ being trivial, assume $d \geq 2$ and that we have already proved the case $d-1$. Then the inductive assumption gives that we have already extend our embedding  $ \varphi :  S_{d-1}(\mathcal{E}) \to S\Xi^{(d-2)}_{[\alpha]}\otimes V$. Then we shall now make an induction on the rank of $\mathcal{E}/S_{d-1}(\mathcal{E})$. \\
   First assume that this rank is $1$. Then let $e \in \mathcal{E}$ which is send to the the standard generator of $\mathcal{E}/S_{d-1}(\mathcal{E}) \simeq \mathcal{E}_\beta$ (so $(a -\beta b)e$ is in $S_{d-1}(\mathcal{E})$) where $\beta$ is in $\alpha + \mathbb{N}$. Note that $[\alpha]$ is in $\mathscr{A}$ because $-\mathscr{A}$ contains the class modulo $\mathbb{Z}$ of any
root of the Bernstein polynomial of $\mathcal{E}$. As we assume that $\mathcal{E}$ has simple pole, $(a - \beta b)e$ is in $S_{d-1}(\mathcal{E}) \cap b\mathcal{E} = bS_{d-1}(\mathcal{E})$ since $S_{d-1}(\mathcal{E})$ is normal in $\mathcal{E}$. Then $\varphi(a - \beta b)e)$ is in $bS\Xi^{(d-2)}_{[\alpha]}\otimes V$ and applying Lemma \ref{manquant 1} we may find $\varepsilon \in S\Xi^{(d-1)}_{[\alpha]}\otimes V$ such that $(a - \beta b)\varepsilon = \varphi(e)$. Then we can extend $\varphi$ to an embedding of $\mathcal{E}$ by sending $e$ to $\varepsilon$ as in the proof of the Embedding Theorem.\\
To complete our induction on the rank of $\mathcal{E}/S_{d-1}(\mathcal{E})$, we have to prove the case of  the rank of $\mathcal{E}/S_{d-1}(\mathcal{E})$ is equal to $k \geq 2$, assuming that the case of rank $k-1$ is already proved, \\
Let $\mathcal{F}$ be a co-rank $1$ normal sub-module containing $S_{d-1}(\mathcal{E})$. This is easily obtained by considering a J-H. sequence  for $\mathcal{E}/S_{d-1}(\mathcal{E})$. Then we have $S_{d-1}(\mathcal{F}) = S_{d-1}(\mathcal{E})$ and $S_d(\mathcal{F}) = \mathcal{F}$. The rank of $\mathcal{F}/S_{d-1}(\mathcal{F})$ is $k-1$ so our inductive assumption gives an embedding $\varphi : \mathcal{F} \to S\Xi^{(d-1)}_{[\alpha]} \otimes V$. Define $\gamma \in \mathbb{Q}^{*+}$ by $\mathcal{E}/\mathcal{F} \simeq \mathcal{E}_\gamma$ and as in the proof of the Embedding Theorem (note that $\mathcal{F}$ has a simple pole because it is normal in $\mathcal{E}$ which has a simple pole) let $e \in \mathcal{E}$ inducing the standard generator of $\mathcal{E}_\gamma$ via the quotient map $\mathcal{E} \to \mathcal{E}/\mathcal{F}$. Then $(a - \gamma b)e$ is in $\mathcal{F}$ and in fact in $b\mathcal{F}$ using the simple pole assumption, so that $\varphi(e)$ is in $bS\Xi^{(d-1)}_{[\alpha]}\otimes V$ and, thanks to Lemma \ref{manquant 1} we may find $\varepsilon \in S\Xi^{(d)}_{[\alpha]} \otimes V)$ which satisfies $(a - \gamma b)\varepsilon = \varphi(e)$. Then, as in the proof of the Embedding Theorem, this allows to define an extension $\tilde{\varphi} : \mathcal{E} \to S\Xi^{(d)}_{[\alpha]}\otimes V$ by putting $\tilde{\varphi}(e) = \varepsilon$. This extension is injective because its kernel has rank at most $1$ so is contained in $S_1(\mathcal{E}) \subset \mathcal{F}$. Moreover, if $\varepsilon$ is not in $S\Xi^{(d-1)}_\alpha\otimes V$ then $B[a].\varepsilon$ is a rank $d+1$ theme, thanks to Lemma \ref{theme 3}, and then $\mathcal{T} := B[a].e \subset \mathcal{E}$ is a $d+1$-theme in $\mathcal{E}$, so $d(\mathcal{T}) = d+1$  contradicting the fact that $S_d(\mathcal{E}) = \mathcal{E}$. So $\tilde{\varphi}$ is an embedding of $\mathcal{E}$ in $S\Xi^{(d-1)}_\alpha\otimes V$ concluding the proof.$\hfill \blacksquare$\\

   \section{Monodromy and the semi-simple filtration}
  
  \subsection{Monodromy}
  
  The goal of the present sub-section is to define the action of the nilpotent part of the logarithm of the  monodromy (logarithm given by $2i\pi b^{-1}a$)  on a simple pole geometric (a,b)-module 
  and to show that the semi-simple filtration of any geometric (a,b)-module $\mathcal{E}$ coincides with the filtration induced on $\mathcal{E}$ by the successive kernels of this nilpotent part  acting on $\mathcal{E}^\sharp$.\\

The first remark is that in a simple pole (a,b)-module the $\mathbb{C}$-linear (bijective) endomorphism $u := b^{-1}a$ satisfies the following commutations relations:
$$ ua - au = a \quad {\rm and} \quad ub - bu = b .$$

\begin{lemma}\label{fond. 1}
Let $X$ be a $\mathbb{C}$-algebra with unit and let $u$ and $x$ elements in $X$ satisfying $ ux - xu = x $. Then for each $n \in \mathbb{N}$ we have $u^nx - xu^n = x((1+u)^{n} - u^{n})$.
\end{lemma} 

\parag{Proof} The result is clear for $n=0, 1$ so assume it is already proved for $n-1$ with $n \geq 2$. Then we have:
\begin{align*}
& u^{n}x  -xu^n = u(xu^{n-1} + x(1+u)^{n-1} - xu^{n-1}) - xu^n  \\
&  u^{n}x  -xu^n  = ux(1+u)^{n-1} - xu^{n} =  x(1+u)^n  - xu^n \qquad  \qquad \qquad \qquad \qquad \qquad \qquad  \blacksquare \\
\end{align*}

Note that if $z$ is a complex number then, replacing $x$ by 
$ z^nx$ gives 
 $$(zu)^nx - x(zu)^n = x((z(u+1))^{n} - (zu)^n) \quad \forall n \in \mathbb{N} .$$ 

Assume now that the series $\sum_{j=0}^\infty (2i\pi u)^n/n! $ converges in the algebra $X$ having a topology for which the product is continuous. Then we have
\begin{align*}
&  \exp(2i\pi u)x - x\exp(2i\pi u) = \sum_{j=0}^\infty [(2i\pi u)^n/n!]x - x[(2i\pi u)^n/n!] \\
& \qquad  = x\big[\sum_{j=0}^\infty (2i\pi (u+1))^n/n! -\sum_{j=0}^\infty (2i\pi u)^n/n!\big] \\
& \qquad  = x\big[\exp(2i\pi(u+1)) - \exp(2i\pi u)\big] = 0,
\end{align*}
which gives  that $\exp(2i\pi u)$ and $x$ commute. 

\begin{lemma}\label{action T 0}
Let $\alpha$ be a rational number in $]0,1]$ and $N$ a non negative integer. Then the $\mathbb{C}$-linear map $\mathcal{T} :=\exp(2i\pi b^{-1}a) := \sum_{q=0}^\infty (2i\pi b^{-1}a)^q/q!$ is well defined on $\Xi_{[\alpha]}^N$ and on  $S\Xi_{[\alpha]}^N$ and is $B[a]$-linear (and also $A$-linear). It is induced by the standard monodromy around $0$  given by $Log\, s \mapsto Log\, s + 2i\pi$. 
\end{lemma}

\parag{Proof} Define $\mathcal{N}$ and $\Delta$ as $\mathbb{C}$-linear endomorphisms of $S\Xi_{\alpha}^{(N)}$ by the formulas  $b^{-1}a = (\Delta + \mathscr{N})$ and 
 $$\Delta(s^{\beta-1}\frac{(Log \, s)^j}{j!}) = \beta s^{\beta-1}\frac{(Log \, s)^j}{j!}$$
 for any $\beta \in \alpha + \mathbb{N}$. Then we have we have  $\mathcal{N}^{N+1} = 0$ since we have for $j \geq 1$
 $$ [\Delta, b^{-1}a](s^{\beta-1}\frac{(Log \, s)^j}{j!}) = -(\beta - 1)s^{\beta-1}\frac{(Log \, s)^{j-1}}{(j-1)!}$$
 and $[\Delta, b^{-1}a](s^{\beta-1}) = 0 $ for $\alpha \not= 1$ and $[\Delta, b^{-1}a](s^{\beta-1}(Log \, s)) = -(\beta - 1)s^{\beta-1} = 0$ for $\alpha = 1$.\\
 Note that $\mathcal{N}(s^{\beta-1}\frac{(Log \, s)^j}{j!}) = s^{\beta-1}\frac{(Log \, s)^{j-1}}{(j-1)!}$, so it commutes with $\Delta$. It also commutes with $a$ and $b$ , then it and is $B[a]$-linear, since we have  $[\Delta, a] = a = [b^{-1}a, a] $ and $[\Delta, b] = b = [b^{-1}a, b]$.\\
  Now   $\exp(2i\pi\Delta) = \exp(2i\pi \alpha)$ on $S\Xi_{[\alpha]}^N$ and  we obtain:
$$ \exp(-2i\pi \alpha) \exp(2i\pi b^{-1}a) = \exp(2i\pi  \mathscr{N}) = \sum_{p=0}^{N+1} \frac{(2i\pi  \mathscr{N})^p}{p!} .$$
Evaluation at  $s^{\beta-1} \frac{(Log \, s)^j}{j!}$ gives, since by definition $\mathcal{T} := \exp(2i\pi b^{-1}a)$:
$$\mathcal{T}(s^{\beta-1}  \frac{(Log \, s)^j}{j!}) = \exp(2i\pi\alpha) \sum_{p=0}^{j} (2i\pi)^p\  s^{\beta-1}  \frac{(Log \, s)^{j-p}}{(j-p)!p!} = \exp(2i\pi \alpha)s^{\beta-1}  \frac{((Log\, s)+ 2i\pi)^j}{j!} $$
thanks to the binomial formula. \\
Then we obtain, since $\exp(2i\pi\mathscr{N}) = (\exp(-2i\pi \alpha))\mathcal{T}$, the equality:
$$\exp(2i\pi\mathscr{N})\big(s^{\beta-1}  \frac{(Log \, s)^j}{j!}\big) =  s^{\beta-1}  \frac{(Log \, s + 2i\pi)^j}{j!}. \qquad \qquad \qquad \qquad \qquad \qquad \blacksquare$$

\bigskip

\begin{cor}\label{ action T 1} For any simple pole geometric (a,b)-module the automorphism $\mathcal{T} := \exp(2i\pi b^{-1}a)$ is well defined and is $B[a]$-linear (and $B[a]$-linear). It is compatible with any $B[a]$-linear map between geometric simple pole (a,b)-modules.
\end{cor}

\parag{Proof} The preceding lemma extends immediately to any $S\Xi_{\mathscr{A}}^{(N)} \otimes V$ and then, using the Embedding Theorem \ref{ET} the result is clear.$\hfill \blacksquare$\\

\begin{defn}\label{fond. 2}
The $B[a]$-linear automorphism $\mathcal{T} := \exp(2i\pi b^{-1}a)$ of a simple pole geometric (a,b)-module $\mathcal{E}$ is called the {\bf monodromy automorphism} of $\mathcal{E}$. 
\end{defn}

\subsection{A direct construction in the formal case}

We give now a direct approach to the monodromy of a formal  simple pole geometric (a,b)-module which does not use the Embedding Theorem. The convergent case seems more difficult to obtain in this way because the convergence of the series defining $\exp(2i\pi b^{-1}a)$ for a simple pole {\it convergent} (a,b)-module is not obvious.
Recall that we denote  $\A := \hat{B}[a]$.

\begin{lemma}\label{mono-1}
Let $E$ be a simple pole $[\alpha]$-primitive  formal   (a,b)-module which is geometric, where $\alpha$ is in $\mathbb{Q} \cap]0,1]$. Define $N := \exp(-2i\pi(\alpha -b^{-1}a)) -1$. Then for each $x \in E$ 
the series
$$2i\pi\mathcal{N}(x) :=  \sum_{p=1}^\infty  (-1)^{p} \frac{N^p}{p}(x) $$
converges in $E$ and $\mathcal{N} : E \to E$ is a $\A$-linear endomorphism of $E$  which satisfies $\mathcal{N}^k = 0$ where $k$ is the rank of $E$.
\end{lemma}

\parag{proof} Since $\exp(-2i\pi(\alpha -b^{-1}a))$ is a unipotent automorphism of $E/bE$, $N$ is nilpotent on $E/bE$ and then $N^k(E) \subset bE$ where $k$ is the rank of $E$. So the series converges for the $b$-adic filtration and $\mathcal{N}(x)$ is well defined for any $x \in E$. The commutation relations $[b^{-1}a, a] = a$ and  $ [b^{-1}a, b] = b$ in the $\mathbb{C}$-linear algebra endomorphisms of $E$ implies the $\A$-linearity of $N$ and then, of $\mathcal{N}$.\\
Then, since $\mathcal{N}$ is nilpotent on $E/bE$ the $b$-linearity implies that $\mathcal{N}^k = 0$ on $E$. Moreover if $f : E \to F$ is a $\A$-linear map between two geometric simple poles formal  (a,b)-modules which are $\alpha$-primitive, the fact that 
$$bf(b^{-1}ax) = f(ax) = af(x) = b(b^{-1}a)f(x)$$
 implies that $f(\mathcal{N}(x)) = \mathcal{N}(f(x))$. 

\begin{defn}\label{mono-2}
Let $E$ be a simple pole geometric formal  (a,b)-module. We define the nilpotent endomorphism $\mathcal{N}$ on $E$ using the direct sum decomposition
$$ E = \oplus_{\alpha \in \mathscr{A}} E_{[\alpha]}$$
where $\mathscr{A}$ is the subset of $\mathbb{Q} \cap ]0,1]$ of class modulo $\mathbb{Z}$ of the eigenvalues of the action of $b^{-1}a$ on $E/bE$. Then $\mathcal{N}$ is the direct sum of the various $\mathcal{N}$ for each $E_{[\alpha]}$. It satisfies $\mathcal{N}^k = 0$ where $k$ is the rank of $E$.  \\
Then we define the $\mathbb{C}$-linear automorphism $\Delta$ of $E$ by the formula
$$ \Delta := b^{-1}a - \mathcal{N}.$$
\end{defn}

\parag{Remarks}
\begin{enumerate}
\item Since $b^{-1}a$ and $\mathcal{N}$ commute, $\Delta$ commutes with $b^{-1}a$ and  $\mathcal{N}$.
\item The commutation relations $[b^{-1}a, a] = a$, $ [b^{-1}a, b] = b$ and the $\A$-linearity of $\mathcal{N}$ implies the commutation relations
$$ [\Delta, a] = a \quad {\rm and} \quad [\Delta, b] = b .$$
\item For any $\A$-linear map $f : E \to F$ between simple poles geometric (a,b)-modules we have $f\circ\Delta = \Delta \circ f$. 
\item $\Delta$ is bijective as a consequence of the fact that $b^{-1}a$ is bijective and the nilpotency of $\mathcal{N}$ proved in the lemma above.
\end{enumerate}

\parag{Conclusion} Let $E$ a simple pole geometric formal (a,b)-module. Defining the monodromy automorphism on $E_{[\alpha]}$, the $[\alpha]$-primitive part of $E$ by the formula
$$ \exp(-2i\pi \alpha)\mathcal{T} = \exp(2i\pi \mathcal{N})$$
the Decomposition Theorem allows to define $\mathcal{T}$ on $E$ as the direct sum of the monodromies of the $E_{[\alpha]}$ where $\alpha$ describes the image of the opposite of the roots of the Bernstein polynomial of $E$ in $\mathbb{Q}/\mathbb{Z}$.\\

\subsection{Nilpotent order}

We come back to the study of {\it convergent}  geometric (a,b)-modules.

\begin{lemma}\label{fond. 3}
Let $\mathcal{E}$ be a simple pole geometric (a,b)-module which is $[\alpha]$-primitive. Then $\mathscr{N}^d = 0$ if $d = d(\mathcal{E})$ is the nilpotent order of $\mathcal{E}$. Moreover $\mathscr{N}^{d-1} \not= 0$.
\end{lemma}

\parag{Proof} First of all remark that $\mathcal{T}$ is an automorphism of $\mathcal{E}$ so that it sends $S_j(\mathcal{E})$ bijectively on itself for any $j$. Moreover, for any normal sub-module $\mathcal{F}$
of $\mathcal{E}$ (so, in particular,  for any $S_j(\mathcal{E})$) $\mathcal{T}$ induces the monodromy automorphism of $\mathcal{F}$ because $\mathcal{F}$ is stable by $b^{-1}a$.\\
Now when $\mathcal{E}$ is semi-simple, as the monodromy of any  rank $1$ simple pole (a,b)-module $\mathcal{E}_\beta$ for $\beta \in \alpha + \mathbb{Z}$ is the product by $\exp(2i\pi\alpha)$, $\mathcal{T}$ is also the mutiplication by $\exp(2i\pi\alpha)$ and the $\mathscr{N} = 0$. \\
When $\mathcal{E}$ is any simple pole geometric (a,b)-module which is $[\alpha]$-primitive, any quotient $S_j(\mathcal{E})/S_{j-1}(\mathcal{E})$ is also $[\alpha]$-primitive  and this implies that $\mathscr{N}(S_j(\mathcal{E})) \subset S_{j-1}(\mathcal{E})$ for any $j$ and then $\mathscr{N}^d = 0$ for $d = d(\mathcal{E})$.\\
To prove the second part of the lemma, remark that the endomorphism $\mathscr{N}$ of $S\Xi^{(N)}_{[\alpha]}\otimes V$ has its image in $S\Xi^{(N-1)}_{[\alpha]}\otimes V$. So, as an $B[a]$-linear embedding of a simple pole (a,b)-module $\mathcal{E}$ (then necessarily geometric and $[\alpha]$-primitive)  in  $S\Xi^{(N)}_{[\alpha]}\otimes V$ commutes with the respective monodromies, using the case $N = d(\mathcal{E})-1$ thanks to Proposition \ref{Manquant 2}, we see that $\mathscr{N}^{d-1}(\mathcal{E}) = 0$ implies that $\mathcal{E} \subset S\Xi^{(d-2)}_{[\alpha]}\otimes V$ which forces  $d(\mathcal{E}) \leq d-1$, contradicting our assumption that $d(\mathcal{E}) = d$. This conclude the proof.$\hfill \blacksquare$\\

Let $\mathcal{E}$ be a simple pole geometric (a,b)-module. Then there exists a finite subset $\mathscr{A}$ in $\mathbb{Q} \cap ]0,1]$ such that $\mathcal{E}$ is  $\mathscr{A}$-primitive. Then define the $B[a]$-linear endomorphism $\mathcal{N}$ of $\mathcal{E}$ by using the direct sum decomposition of $\mathcal{E}$ (see Theorem \ref{DT}):
$$ \mathcal{E} \simeq \oplus_{\alpha \in \mathscr{A}} \ \mathcal{E}_{[\alpha]} $$
Then define  $\mathscr{N}$ as the direct sum of the endomorphism $\mathscr{N}$ on each $\mathcal{E}_{[\alpha]}, \alpha \in \mathscr{A}$.\\
Define also $\Delta$ on $\mathcal{E}$  as $b^{-1}a + \mathscr{N}$.\\
Then the following theorem is an easy consequence of the previous lemma using the fact that $S_d(S\Xi^{(N)}_{[\alpha]}\otimes V) = \Xi^{(d-1)}_{[\alpha]}\otimes V$ and Lemma \ref{basic filt.ss}.

\begin{thm}\label{fond. 4}
Let $\mathcal{E}$ be a geometric (a,b)-module and let $\mathscr{N}$ be the nilpotent part of the monodromy acting on $\mathcal{E}^\sharp$. Then intersection with $\mathcal{E}$ of  the kernel of $\mathscr{N}^j$ is equal to $S_{j}(\mathcal{E})$ for all $j \in [0,d]$. So $d= d(\mathcal{E}) = d(\mathcal{E}^\sharp) $ is the nilpotent order of the action of the monodromy on $\mathcal{E}^\sharp$. $\hfill \blacksquare$\\
\end{thm}

In other terms the previous corollary explains that the {\bf semi-simple filtration} of a geometric (a,b)-module coincides with the filtration induced on $\mathcal{E}$ by the successive kernels of the powers of  nilpotent part  of the monodromy of $\mathcal{E}^\sharp$.\\
Note that, in general, the inclusion  $\mathscr{N}(S_j(\mathcal{E})) \subset S_{j-1}(\mathcal{E})$  is not an equality for a simple pole geometric (a,b)-module.
For instance if $\mathcal{E} = \mathcal{F} \oplus \mathcal{G}$ where $\mathcal{F}$ is semi-simple and $\mathcal{G}$ is not semi-simple with $d(\mathcal{G}) = 2$, $\mathscr{N}(\mathcal{E}) = \mathscr{N}(\mathcal{G}) \subset S_1(\mathcal{G})$ is strictly contained in $S_1(\mathcal{E}) = \mathcal{F} \oplus S_1(\mathcal{G})$ when $\mathcal{F} \not= \{0\}$.\\

We collect in the following proposition the main tools we  have obtained to compute the nilpotent order of a geometric (a,b)-module.(see \cite{[B.09]} for the formal case; the convergence case is analogous).
 \begin{prop}\label{caract. ss}
    Let \ $\mathcal{E}$ \ be a geometric (a,b)-module. Then we have the following properties :
  \begin{enumerate}[i)]
  \item For each subset  $\mathscr{A}  \in \mathbb{Q}\cap ]0,1]$ we have $S_j(\mathcal{E}_{[\mathscr{A}]}) = S_j(\mathcal{E})_{[\mathscr{A}]}$, for each $j \geq 1$.
  \item Any \ $[\alpha]-$primitive  sub-theme \ $\mathcal{T}$ \ in \ $\mathcal{E}$ \ of rank \ $j$ \ is contained in \ $S_j(\mathcal{E})$.
  \item Any  \ $[\alpha]-$primitive quotient theme \ $\mathcal{T}$ \ of \ $S_j(\mathcal{E})$ \ has  rank \ $\leq j$. 
  \item The nilpotent order of \ $\mathcal{E} = \mathcal{E}_{[\mathscr{A}]}$ \ is equal to \ $d$ \ if and only if \ $d$ \ is the maximal rank of an \ $[\alpha]-$primitive quotient theme of \ $\mathcal{E}$ where $\alpha$ is in $\mathscr{A}$.
  \item The nilpotent order  of \ $\mathcal{E} = \mathcal{E}_{[\mathscr{A}]}$ \ is equal to \ $d$ \ if and only if \ $d$ \ is the maximal rank of an \ $[\alpha]-$primitive sub-theme of \ $\mathcal{E}$ where $\alpha$ is in $\mathscr{A}$.$\hfill \blacksquare$\\
  \end{enumerate} 
  \end{prop}
  
   \begin{lemma}\label{nilpotent 1}
Let $\mathcal{F}$ and $\mathcal{G}$ be two sub-modules of a  geometric (a,b)-module. Assume that $d(\mathcal{F}) \leq p_1$ and $d(\mathcal{G}) \leq p_2$. Then $d(\mathcal{F} + \mathcal{G}) \leq p$ where $p = \sup \{p_1, p_2\}$.
\end{lemma}

\parag{Proof} Assume that we have a surjective (a,b)-linear map $\pi : \mathcal{F} + \mathcal{G} \to T$ where $T$ is an $[\alpha]$-primitive theme of rank $q > p$. Then let $e$ be a generator of $T$ as a $B[a]$-module and let $u \in \mathcal{F}$ and $v \in \mathcal{G}$ such that $e = \pi(u + v)$. Note $T_1 := \pi(B[a] u)$ and $T_2 := \pi(B[a] v)$. Let $\tilde{T}_i$ be the normalization of $T_i$ in $T$ for $i = 1, 2$. Then if $\tilde{T}_i \not= T$ for $i = 1, 2$ these two normal sub-module are contained in the co-rank $1$ sub-module of $T$, and this is not possible since $T_1 + T_2 = T$.\\
So we have, for instance $\tilde{T}_1 = T$. But this means that there exists $n \in \mathbb{N}$ such that $b^ne$ is in $T_1$. Then $b^nT$ is contained in $T_1$ and so $\mathcal{F}$ contains the sub-module $\pi^{-1}(b^nT)$ which admits $b^nT$ as a quotient. But $b^nT$ is a rank $q$ $[\alpha]$-primitive theme and then $d(B[a] u) \geq q$. This contradicts our hypothesis that $d(\mathcal{F}) \leq p$ since we assume that $q > p$.\\
So any quotient theme of $\mathcal{F} + \mathcal{G}$ has rank at most equal to $p$, concluding the proof, thanks to Proposition \ref{caract. ss}.$\hfill \blacksquare$\\

\begin{cor}\label{nilpotent 2}
Let $\mathcal{E}$ be a geometric (a,b)-module. For each integer  $j \geq 0$, $S_j(\mathcal{E})$ is the subset of all $x \in \mathcal{E}$ such that $d(B[a] x) \leq j$.
\end{cor}

\parag{Proof} It is clear that $x \in S_j(\mathcal{E})$ implies $d(\Abc x) \leq j$. Conversely, let $x \in \mathcal{E}$ with $d(\Abc x) \leq j$. Then, thanks to the previous lemma,  we have
$$d(\Abc x + S_j(\mathcal{E})) \leq j$$
 and this implies that $\Abc x + S_j(\mathcal{E}) \subset  S_j(\mathcal{E})$ because $ S_j(\mathcal{E})$ is the maximal sub-module in $\mathcal{E}$ with nilpotent order equal to $j$ (see Remark 2. following Lemma \ref{basic filt.ss}). This is enough to conclude.$\hfill \blacksquare$\\

\end{document}